\documentclass[conference]{IEEEtran}
\IEEEoverridecommandlockouts
\usepackage{graphicx} 
\usepackage{amsmath}
\usepackage{amsthm}
\usepackage{amssymb}
\usepackage{xcolor}
\usepackage{verbatim}
\usepackage[hidelinks]{hyperref}
\usepackage{algorithm}
\usepackage[noend]{algpseudocode}
\usepackage[bottom=2.0 cm,top=2.0 cm, left=2.0 cm, right=2.0 cm]{geometry}
\usepackage{bbm}
\usepackage[bb=boondox]{mathalfa}
\usepackage{graphicx} %
\hypersetup{
    colorlinks=true,
    allcolors=black 
}

\newcommand{\X}{\mathbf{X}}

\newcommand{\Y}{\mathbf{Y}}
\newcommand{\D}{\mathbf{D}}

\newcommand{\PL}{\text{\L}}

\definecolor{xkchen}{rgb}{0.75, 0.2, 1.0}
\newtheorem{assumption}{Assumption}

\newtheorem{theorem}{Theorem}
\newtheorem{corollary}{Corollary}[theorem]
\newtheorem{lemma}{Lemma}
\newtheorem{remark}{Remark}

\usepackage{cite}
\usepackage{amsmath,amssymb,amsfonts}
\usepackage{graphicx}
\usepackage{textcomp}
\usepackage{xcolor}
\def\BibTeX{{\rm B\kern-.05em{\sc i\kern-.025em b}\kern-.08em
    T\kern-.1667em\lower.7ex\hbox{E}\kern-.125emX}}
\begin{document}

\title{\vspace{-0.7cm}A Line-search-free Method for Adaptive Decentralized Optimization
\thanks{Chen and Scutari are with the School of Industrial Engineering, Purdue University; emails:\texttt{<chen4373,gscutari>@purdue.edu}. }
}
\author{Xiaokai Chen, Ilya Kuruzov, and Gesualdo Scutari}\vspace{-1.8cm}

\maketitle

\begin{abstract}
\label{sec:abstract}
We study decentralized optimization over networks where agents cooperatively minimize a smooth (strongly) convex sum of local losses while communicating only with immediate neighbors. Prevailing decentralized methods require either centralized knowledge of global problem and network parameters for stepsize tuning--hence impractical, or 
costly per-iteration line-searches that demand access to local function values. We propose line-search-free, fully decentralized algorithms in which each agent adapts its stepsize using only past local iterates and gradients--with no extra function evaluations and no global tuning.
  The key technical ingredient is a new Lyapunov function, 
  from which a natural adaptive stepsize rule emerges: at each iteration, each agent selects the largest stepsize that guarantees descent, based solely on a local curvature estimate built from successive gradients.    The proposed algorithms enjoy strong theoretical  guarantees: sublinear convergence rates for merely convex objectives and linear rates under strong convexity. Numerical experiments on standard benchmarks show consistent improvements over the state of the art, both adaptive and non-adaptive.
\end{abstract}


\section{Introduction}
\label{sec:introduction}
We consider  the following convex optimization  problem 
\begin{equation}
    \label{eq:problem}
    \tag{P}
    \begin{aligned}    
&\min_{x\in\mathbb{R}^d}f(x),\quad f(x):=\frac{1}{m}\sum_{i=1}^mf_i(x),
    \end{aligned}
\end{equation}
over a network of $m$ agents, where $f_i$ is the cost function accessible only to agent $i\in[m]:=\{1,\cdots, m\}$, assumed to be (strongly) convex and \textit{locally} smooth. All   agents are embedded in a   communication network  modeled as a fixed, undirected and connected graph.

Decentralized optimization in the form~\eqref{eq:problem} arises   in several scientific and engineering areas, including distributed machine learning, 
multiple-agent control and coordination, 
and sensor network information processing. 
A large body of methods exists for \eqref{eq:problem}; representative methods include gradient tracking~\cite{di2016next,sun2022distributed,nedic2017achieving} and primal–dual schemes~\cite{shi2015extra,yuan2018exact}; see~\cite{nedic2018network,sayed2014adaptation} for  some overviews.  
As in centralized optimization, convergence of decentralized methods hinges on an appropriate stepsize policy--typically a {\it single} stepsize shared across agents--with either  {\bf constant} or {\bf  diminishing} rules.

Under {\bf constant}  stepsize rules, existing theory imposes upper bounds on the stepsize  values that depend on global problem/network constants--e.g., Lipschitz constants of $\nabla f_i$, the spectral gap of the gossip matrix, and related topological descriptors. These quantities are rarely known a priori in practice. Consequently, despite theoretical guidance, practitioners often fall back on manual tuning of a single, common stepsize for each problem/network instance, which can be slow and brittle because it ignores the heterogeneous   geometry of the loss landscape   along the trajectory. 
  {\bf Diminishing stepsizes} have been adopted to partially sidestep unavailable global problem/network-dependent constants, e.g.., ~\cite{chen2012fast,nedic2014distributed,choi2022convergence}. This, however, comes at   the cost of  significantly degraded rates (e.g., sublinear iteration complexity even under strong convexity) and   sensitivity to schedule design. 

 This has recently motivated growing interest in {\it adaptive}  methods that adjust stepsizes based on     loss curvature information while still guaranteeing global convergence. Below we summarize the most relevant literature on the subject.  \vspace{-0.2cm}
\subsection{Literature review}
   The literature for adaptive {\bf centralized} optimization has grown rapidly in the past few years. Several stepsize selection rules based on local estimates of the loss curvature have been proposed, ranging from the classical line-search~\cite{bertsekas1997nonlinear}, Barzilai-Borwein~\cite{barzilai1988two,zhou2025adabb}, and Polyak~\cite{polyak1969minimization} rules, to more recent schemes leveraging past gradient and iterate information~\cite{malitsky2019adaptive,malitsky2024adaptive,malitsky2018first}, as well as methods designed for stochastic and large-scale settings such as AdaGrad~\cite{duchi2011adaptive}, Adam~\cite{kingma2014adam}, AMSGrad~\cite{reddi2019convergence}, and their variants~\cite{li2019convergence,ward2020adagrad}. 
   However, all these methods 
  cannot be implemented over mesh networks without a centralized node.

   In contrast, the literature for adaptive {\bf decentralized} algorithms is fairly scarce.
   Direct adaptations of centralized methods~\cite{malitsky2018first,latafat2024adaptive} to decentralized settings exhibit significant practical limitations:  the stepsize selection either leads to extensive        communications (to verify line-search based terminations) per iteration   or requires global scalar information  not available at the agents' side (see~\cite{chen2025parameter} for details).    
   Other approaches~\cite{atenas2024distributed,aldana2024towards}, despite being adaptive, still require {\it problem/network-dependent} parameter knowledge for stepsize selection.    In~\cite{Chen2026AdaptivePS}, a decentralized Polyak stepsize rule is proposed; however,   its $\mathcal{O}(1/\sqrt{mk})$ ($k$ is the iteration index) rate matches that of non-adaptive methods, so the adaptivity yields no rate improvement, and slow convergence in the strongly convex case. Furthermore, it requires agents'  {losses to be {\it globally Lipschitz} continuous (i.e., $\nabla f_i$ uniformly bounded), which excludes a vast class of problems}. The proposals~\cite{nazari2022dadam,chen2023convergence,li2024problem} focus on the stochastic instance of problem \eqref{eq:problem}; they explore  adaptivity through stochastic gradient normalization, achieving convergence rate of   $\mathcal{O}(1/\sqrt{k})$, which are quite unfavorable in determinist setting considered in this paper.  
   
   Advancements have been made in~\cite{kuruzov2025adaptive,kuruzovachieving,chen2025parameter,chen2026adaptive,Xu2025AnAP}. They proposed decentralized algorithms that adaptively select stepsizes via local backtracking line-searches combined with global/local min-consensus protocols, achieving robust and fast convergence guarantees. These represent the current state of the art;  however, the line-search procedures introduce extra per-iteration function evaluations--hence computational complexity--and the algorithms are not implementable when objective function values are unavailable.

   To summarize, no existing adaptive decentralized method  provides robust convergence guarantees while avoiding the computational overhead of line-searches. This paper addresses this gap;  our main contributions are as follows. \vspace{-0.1cm}

\vspace{-1mm}
   \subsection{Main contributions}
\label{subsec:contribution}
\textit{1. Algorithm design: Line-search-free adaptive algorithms.}  
We propose two   decentralized algorithms (Algorithm~\ref{alg:line-search-free} and~\ref{alg:line-search-free_local}) for~\eqref{eq:problem} that adapt each agent's stepsize without line-searches and without knowledge of any global problem or network constant. Algorithm~\ref{alg:line-search-free} relies on a lightweight global scalar averaging protocol, while Algorithm~\ref{alg:line-search-free_local} requires only local (neighbor-to-neighbor) scalar exchanges, making it implementable on any connected network. 

 \textit{2. Key design principle.} The algorithms are derived from a strengthened convergence analysis of the 
Condat-V\~{u} splitting~\cite{condat2013primal,vu2013variable}, which is of independent interest.  We propose a novel   augmented Lyapunov 
function that {\it incorporates trajectory information}--consecutive iterate 
differences and primal Lagrangian gaps. Unlike standard Lyapunov functions used in the literature, whose descent conditions couple the stepsize to both smoothness and network spectral quantities--not computable by the agents, ours yields a descent condition that depends only on local curvature estimates built from past gradients and iterates. Each agent can then maximize  its local stepsize subject to this 
computable descent condition, obtaining adaptivity as a byproduct of the convergence mechanism itself rather than as a separate heuristic layer.  


\textit{3. Strong convergence guarantees:} 
For merely convex losses $f_i$, the algorithms  are proved to reach  an $\varepsilon$-solution (measured by a suitable metric) of Problem~\eqref{eq:problem} in $\mathcal{O}(\widetilde{L}/\varepsilon)$  iterations (communications), where $\widetilde{L}$  is the Lipschitz constant of  $\nabla f_i$'s {\it restricted to   the   convex hull} of the trajectory generated by the algorithm.  For strongly convex $f_i$'s,   
  linear convergence  is certified with communication complexity in the order of $\mathcal{O}\left(\frac{\widetilde{\kappa}}{1-\lambda_2(W)}\log(1/\varepsilon)\right)$, where $\lambda_2(W)$ is the second largest eigenvalue of  the gossip matrix $W$   and $\widetilde{\kappa}$ is the   condition number of $f_i$'s {\it restricted to  the  convex hull} of the iterates.    Quite remarkably, both rates outperform those of existing nonadaptive algorithms, since  $\widetilde{L}$   and    $\widetilde{\kappa}$  can be significantly smaller than their global counterparts. Our numerical results validate these theoretical findings.

\section{Algorithm Design and Convergence Analysis}
\label{sec:problem}
We study~\eqref{eq:problem} under the following standard assumptions.
\begin{assumption}[Objective function]
\label{ass:function}
    Each   function $f_i:\mathbb{R}^d\rightarrow \mathbb{R}$ is continuously differentiable,   \textit{locally} smooth, and   $\mu$-strongly convex, with $\mu\in[0,\infty)$. If $\mu=0$,   $f$   admits a minimizer over $\mathbb R^d$.
    \end{assumption}

Unlike most of     decentralized optimization literature (see, e.g., ~\cite{nedic2018network,sayed2014adaptation} and references therein), this work  does {\it not}require {\it global} smoothness of the agents' losses.  This permits covering  a broader class of practical problems, including  Poisson linear inverse problem~\cite{bauschke2017descent} and  Poisson Regression~\cite{atenas2025distributed}, to name a few.  
 
\begin{assumption}[Communication graph] \label{ass:grap} The   graph  $\mathcal G=(\mathcal{V},\mathcal{E})$ is undirected and connected, where   $\mathcal{V}:=[m]:=\{1,\ldots, m\}$ is the set of agents and $(i,j)\in\mathcal{E}$ if and only if there is a communication link between agents $i$ and $j$.
\end{assumption}
 We denote by  $\mathbb{S}^m_{+}$ (resp. $\mathbb{S}^m_{++}$)   the set of $m\times m$  real  symmetric positive semidefinite  (resp. definite) matrices.

We begin   rewriting~\eqref{eq:problem} in an equivalent lifted form. Introducing local copies $x_i \in \mathbb{R}^d$ of the  variable  $x\in \mathbb{R}^d$---the $i$-th controlled by agent $i$---the row-wise stack matrix $\mathbf{X}:=[x_1,x_2,\ldots,x_m]^\top\in\mathbb{R}^{m\times d}$, and the    lifted  functions $$F(\X):=\sum_{i=1}^m f_i(x_i), \quad  G(\X):=\delta_{\{0\}}(\X),$$
where $\delta_{\{0\}}:\mathbb{R}^{m\times d}\rightarrow\mathbb{R}\cup\{\infty\}$ is the indicator function of $\{0\}$, Problem~\eqref{eq:problem} can be equivalently written as:\vspace{-0.1cm}
\begin{equation}
    \label{eq:primal}
    \tag{P$^\prime$}
   \min_{\X\in\mathbb{R}^{m\times d}} F(\X)+G(\PL\X),
\end{equation}
where $\PL\in\mathbb{S}_+^m$ is an $m \times m$ matrix satisfying $\texttt{null}(\PL)=\texttt{span}(\mathbf{1}_m)$. This condition ensures  that any  $\X$ satisfying $\PL\X=\mathbf{0}$ is consensual, i.e., $\X=\mathbf{1}_m x^\top$ for some  $x\in\mathbb{R}^d$.
A standard choice for $\PL$ in  decentralized algorithms is $\PL=(I-\widetilde{W})^{1/2}$ for some $\widetilde{W}\in\mathcal{W}_{\mathcal{G}}$, where $\mathcal{W}_{\mathcal{G}}$ is the class of gossip matrices defined next.
\begin{assumption}
    [Gossip matrices]
    \label{ass:W}
    $\mathcal{W}_{\mathcal{G}}$ denotes the set of symmetric, doubly stochastic  matrices $\widetilde{W}:=(\widetilde{w}_{i,j})_{i,j=1}^m$ that are compliant with $\mathcal{G}$, i.e., $\widetilde{w}_{ii}>0$ for all $i\in[m]$; $\widetilde{w}_{ij}>0$ for all $(i,j)\in\mathcal{E}$; and $\widetilde{w}_{ij}=0$ otherwise.
    \end{assumption}
Several instances of such matrices have been employed in practice; see, e.g.,~\cite{nedic2020distributed,nedic2018network,sayed2014adaptation}. Throughout  the paper  we number the eigenvalues    of matrix $W \in \mathbb{S}^m$ in nonincreasing  order, i.e.,   $\lambda_1(W)\geq\cdots\geq\lambda_m(W)$. 
\subsection{Revisiting convergence of the Condat-V\~{u} splitting method} \label{subsec:revisit}
We start from the Condat-V\~{u} splitting~\cite{condat2013primal,vu2013variable}---a classical primal-dual method for  problems in the form of~\eqref{eq:primal}.   Given parameters   $\alpha,\gamma,\sigma>0$ (to be properly determined) and an initialization $\X^0,\X^{-1},\Y^0\in\mathbb{R}^{m\times d}$, the algorithm reads: for $k=0,1,2,\ldots,$
\begin{equation}
    \label{eq:original condat-vu}
\begin{aligned}
\Y^{k+1}&=\texttt{prox}_{\sigma\alpha G^*}(\Y^k+\sigma\alpha\PL((1+\gamma)\X^k-\gamma\X^{k-1}));\\
    \X^{k+1}&= \X^k -\alpha(\nabla F(\X^k)+\PL\Y^{k+1}),
\end{aligned}
\end{equation}
where $\texttt{prox}_{\sigma\alpha G^*}(\X):=\arg\min_{\Y}G^*(\Y)+\frac{1}{2\sigma\alpha}\|\X-\Y\|^2$ denotes the proximal operator of the conjugate function $G^*$ of $G$ (applied row-wise), scaled by $\sigma\alpha$. 

Using algorithm~\eqref{eq:original condat-vu} with the existing convergence machinery~\cite{ryu2022large} presents  two fundamental obstacles.

\smallskip
 \textit{\textbf{(i)}  The classical stepsize condition is impractical for adaptive methods.}
Setting  $\PL=(I-\widetilde{W})^{1/2}$ for some $\widetilde{W}\in\mathcal{W}_{\mathcal{G}}$,   the classical convergence analysis~\cite{ryu2022large}---which requires the stronger assumption that $F$ is globally $L$-smooth (hence violating Assumption~\ref{ass:function})---guarantees that  $\{\X^k\}_{k\geq 0}$ generated by~\eqref{eq:original condat-vu} converges to a solution $\X^*$ of~\eqref{eq:primal}, provided\vspace{-0.1cm}
\begin{equation}
    \label{eq:original stepsize}
    \alpha L/2+\sigma\alpha^2\|I-\widetilde{W}\|<1. 
\end{equation}
This condition suffers from several limitations. First, it couples the stepsize $\alpha$ and the tuning parameter $\sigma>0$ with the network-dependent quantity $\|I-\widetilde{W}\|$, which is typically unavailable to individual agents. One could  upper bound $\|I-\widetilde{W}\|\leq 2$ and substitute into~\eqref{eq:original stepsize}; however, this yields a highly conservative stepsize  that leads to slow convergence in practice~\cite{chen2025parameter}. Second, the dependence on the \textit{global} smoothness constant $L$ (unknown to the agents) is restrictive, as $L$ can be significantly larger than the \textit{local} smoothness along the trajectory of the algorithm.  Third,~\eqref{eq:original stepsize} involves the auxiliary parameter $\sigma$, whose tuning interacts nontrivially with $\alpha$: making $\alpha$ adaptive while $\sigma$ is fixed (or vice versa) proves problematic, and existing attempts along this route exhibit high sensitivity of the algorithm performance to the choice of $\sigma$~\cite{chen2025parameter}.

\smallskip
\textit{\textbf{(ii)}  Algorithm~\eqref{eq:original condat-vu} is not directly implementable in a decentralized setting.} The updates in~\eqref{eq:original condat-vu} involve the linear operator  $\PL=(I-\widetilde{W})^{1/2}$, whose application to the {$X$-and $Y$-variables}   requires global coordination and cannot be carried out via local communication among neighbors.

\smallskip
Next, we address both limitations above: {\bf (i)} we develop a new convergence analysis for~\eqref{eq:original condat-vu} based on a novel Lyapunov function, yielding a  stepsize condition  that is \textit{independent} of the network parameters (Sec.~\ref{subsec:new analysis}); {\bf (ii)} we  make~\eqref{eq:original condat-vu} implementable in a decentralized fashion and use the network-independent analysis to guide the adaptive selection of the stepsize,  {tailored to the local curvature of the agents' losses, thereby   getting rid of the dependence on the much more conservative global smoothness constant $L$} (Sec.~\ref{subsec:stepsize}).
\vspace{-1mm}
\subsection{Network-independent convergence analysis}
\label{subsec:new analysis}
Let $\mathcal{L}(\X,\Y)$ be the Lagrangian function of Problem~\eqref{eq:primal}:
\begin{equation}
    \label{eq:lagarangian line-search-free}
    \mathcal{L}(\X,\Y)=F(\X)+\langle\PL \X,\Y\rangle-G^*(\Y);
\end{equation}
and, for any given  saddle point $(\X^*,\Y^*)$  of~\eqref{eq:lagarangian line-search-free}, define the     primal gap  at  $\X$ as:
\begin{equation}
    \label{eq:L_p gap}
  \begin{aligned}
      &\Delta  \mathcal{L}_*^k:=\mathcal{L}(\X^k,\Y^*)-\mathcal{L}(\X^*,\Y^*)\\
      =&F(\X^k)-F(\X^*)+\langle \PL\Y^*,\X^k-\X^*\rangle\geq 0,
  \end{aligned}
\end{equation}
where  the inequality follows from the convexity of $F$.

Our enhanced  convergence analysis of algorithm~\eqref{eq:original condat-vu} hinges on  the following novel Lyapunov function evaluated along the iterates $\{(\X^k,\Y^k)\}$ of~\eqref{eq:original condat-vu}: for any given saddle point $(\X^*,\Y^*)$  of~\eqref{eq:lagarangian line-search-free}, let \vspace{-0.2cm}
\begin{equation}\label{eq:Lyap}\begin{aligned}
    &V_*^k:=\|\X^k-\X^*\|^2+\frac{1}{\sigma}\|\Y^k-\Y^*\|^2\\&+\frac{1}{2}\|\X^{k}-\X^{k-1}\|^2+2\gamma\alpha\Delta  \mathcal{L}_\star(\X^{k-1}),
\end{aligned}\end{equation}
where $\alpha,\gamma,\sigma>0$ are parameters to   be properly chosen.  Note that $V^k_\star$ is continuous in its arguments and $V^k_\star\geq 0$, with equality if and only if $(\X^k,\Y^k)=(\X^{k-1},\Y^{k-1})=(\X^*,\Y^*)$;  hence, it is a valid measure of optimality. 

The following result provides conditions for the monotone decrease of $V^k_\star$ along $\{(\X^k,\Y^k)\}$.
\begin{lemma}
\label{lemma:descent lemma condat-vu}
 Let $\{(\X^k,\Y^k)\}$ be the iterates generated by~\eqref{eq:original condat-vu} under Assumptions~\ref{ass:function}-\ref{ass:W}, with the additional requirement that $F$ is $L$-smooth globally, and the following tuning:    $\widetilde{W}\in\mathcal{W}_{\mathcal{G}}\cap\mathbb{S}_{++}^m$, $\sigma>0,$ $\gamma=1$, and \vspace{-0.2cm}\begin{equation}\label{eq:step-size}
     \alpha\leq \frac{1}{\sqrt{L^2+2\sigma}+L}.
 \end{equation}
Then, for any given     saddle point $(\X^*,\Y^*)$ of~\eqref{eq:lagarangian line-search-free}, it holds
\begin{equation}
    \label{eq:descent condat-vu}
   \begin{aligned}
        V_{\star}^{k+1}\leq V_{\star}^k,\quad \forall k=0,1,\ldots .
   \end{aligned}
\end{equation}
\end{lemma}
The requirement that  $\widetilde{W}\succ 0$  can be readily met by using the shifted matrix $W=(1-c)I+c\widetilde{W}\in\mathcal{W}_{\mathcal{G}}$, for some $c\in(0,1/2)$ and   $\widetilde{W}\in\mathcal{W}_{\mathcal{G}}$~\cite{nedic2020distributed,nedic2018network,sayed2014adaptation}.

Crucially, the stepsize condition \eqref{eq:step-size}  depends only on the smoothness constant $L$ and the free parameter $\sigma$, with no appearance of the network quantity $\|I-\widetilde{W}\|$. This decoupling is the key insight that enables the design of adaptive stepsizes based on local curvature alone.
\vspace{-1mm}
\subsection{Decentralization of~\eqref{eq:original condat-vu} and adaptive stepsize selection}
\label{subsec:stepsize}
We now address limitation (ii) of Section~\ref{subsec:revisit}. To make the algorithm~\eqref{eq:original condat-vu}   implementable in a decentralized setting, using $\PL=(I-\widetilde{W})^{1/2}$ with   $\widetilde{W}\in \mathcal{W}_{\mathcal{G}}$,   we introduce the change of  variable $\D^k=\PL\Y^k$,  with   $\D^0=\PL\Y^0=\mathbf{0}$.   Since $G=\delta_{\{0\}}$, we have  $G^*(\Y)=0$ for all $\Y$, and the proximal step on $\Y$ in ~\eqref{eq:original condat-vu}  reduces to a linear update. Using these facts, one can verify  that~\eqref{eq:original condat-vu} is equivalent to  \begin{equation}
    \label{eq:Decentralized condat-vu}
    \begin{aligned}
   \D^{k+1}&=\D^k+\sigma\alpha(I-\widetilde{W})((1+\gamma)\X^k-\gamma\X^{k-1});\\
   \X^{k+1}&=\X^k-\alpha(\nabla F(\X^k)+\D^{k+1}).
    \end{aligned}
\end{equation}
Algorithm~\eqref{eq:Decentralized condat-vu} is now implementable in a decentralized setting; furthermore, unlike  existing adaptive decentralized methods~\cite{kuruzovachieving,kuruzov2025adaptive,chen2025parameter,chen2026adaptive,Xu2025AnAP},  it requires  only one round of communication per iteration.

Building on the network-independent analysis of Lemma~\ref{lemma:descent lemma condat-vu}, we now drop the global $L$-smoothness assumption  and select the stepsize adaptively. We apply a Lyapunov function analogous to $V^k_\star$ to  the decentralized update~\eqref{eq:Decentralized condat-vu}, replacing the fixed parameters $\alpha,\sigma,\gamma$ with non-stationary sequences $\{\alpha^k\}_{k\geq 0}$, $\{\sigma^k\}_{k\geq 0}$, and $\{\gamma^k\}_{k\geq 0}$ that are  adapted along the  iterates. Specifically, given $\X^{-1},\X^0\in\mathbb{R}^{m\times d}$, $\D^0=\mathbf{0}$, and   $\widetilde{W}\in\mathcal{W}_{\mathcal{G}}\cap\mathbb{S}_{++}^m$, we have:  \begin{subequations}
 \label{eq:primal and dual update}
             \begin{align}
\D^{k+1}&=\D^k+\sigma^k\alpha^k(I-W)((1+\gamma^{k})\X^k-\gamma^{k}\X^{k-1});\label{eq:dual update}\\
    \X^{k+1}&=\X^k-\alpha^k(\nabla F(\X^k)+\D^{k+1}),  \label{eq:primal update}           
             \end{align}   
            \end{subequations}
 for $k=0,1,2,\ldots, .$
Recalling that $\D^k=\PL\Y^k$, the update of the auxiliary variable $\{\Y^k\}_{k\geq 0}$ follows from~\eqref{eq:dual update}:  \begin{equation}
    \label{eq:auxiliary update}   \Y^{k+1}=\Y^k+\sigma^k\alpha^k((1+\gamma^{k})\PL\X^k-\gamma^{k}\PL\X^{k-1}),
\end{equation} for  $k=0,1,2,\cdots,$ and  $\Y^0=\mathbf{0}$.

 The  Lyapunov function $V^k_\star$ in~\eqref{eq:Lyap}, adjusted to the non-stationary sequences $\{\alpha^k\}_{k\geq 0}$, $\{\sigma^k\}_{k\geq 0}$, and $\{\gamma^k\}_{k\geq 0}$, becomes (with a slight abuse of notation):   
\begin{equation}
\label{eq:merit function:line-search-free}
\begin{aligned}
        &V^{k}_\star:=\|\X^k-\X^*\|^2+\frac{1}{\sigma^k}\|\Y^k-\Y^*\|^2\\&+\frac{1}{2}\|\X^{k}-\X^{k-1}\|^2+2\gamma^{k}\alpha^{k}\Delta  \mathcal{L}^{k-1}_\star\geq 0,
\end{aligned}
\end{equation} where $\X^k,\X^{k-1}$ and $\Y^k$ are given by (\ref{eq:primal update}) and (\ref{eq:auxiliary update}), respectively. The counterpart of Lemma~\ref{lemma:descent lemma condat-vu} reads as follows.  
\begin{lemma}
    \label{lemma:Lyapunov dynamic}
   Let $\{(\X^k,\D^k)\}$ be the iterates generated by~\eqref{eq:primal and dual update} and $\{\Y^k\}$ the corresponding auxiliary variables defined in~\eqref{eq:auxiliary update}, under Assumptions~\ref{ass:function}-\ref{ass:W} and the following tuning: $W\in\mathcal{W}_{\mathcal{G}}\cap\mathbb{S}_{++}^m$, $\gamma^0=1$, and $\gamma^k=\alpha^k/\alpha^{k-1}$ for  $k\geq 1$. Then, for any  given  saddle point $(\X^*,\Y^*)$ of~\eqref{eq:lagarangian line-search-free} and any $\zeta^k>0$, it holds\vspace{-0.2cm}
    \begin{equation}
    \label{eq:Lyapunov dynamic}
       \begin{aligned}
       &V^{k+1}_\star\\
       \leq& V^k_\star-\left(\frac{1}{2}-\alpha^kL^k-\alpha^k\zeta^k\right)\|\X^{k+1}-\X^{k}\|^2
       \\ &-\left(\frac{1}{2}-\alpha^kL^k\right)\|\X^{k}-\X^{k-1}\|^2\\
       &\!\!\!\!\!-\left(\frac{1}{\sigma^k}-\frac{\alpha^k}{\zeta^k}\right)\|\Y^{k+1}-\Y^k\|^2-\frac{\sigma^{k+1}-\sigma^k}{\sigma^{k}\sigma^{k+1}}\|\Y^{k+1}-\Y^*\|^2
       \\
       &\!\!\!\!-[(2+2\gamma^{k})\alpha^k-2\gamma^{k+1}\alpha^{k+1}]\Delta  \mathcal{L}^k_\star, \quad \forall k\geq 0,
       \end{aligned}
    \end{equation}
    where $L^k$ is a proxy for the local Lipschitz constant of $\nabla f_i$'s, defined as\vspace{-0.3cm}
    \begin{equation}
    \label{eq:L^k}
    L^k:=\frac{\|\nabla F(\X)-\nabla F(\X^{k-1})\|}{\|\X^k-\X^{k-1}\|}.
    \end{equation}
\end{lemma}
\begin{proof}
    See Appendix.~\ref{appx:proof of lemma lyapunov}
\end{proof}

To guarantee descent of the Lyapunov function, i.e., $V^{k+1}_\star\leq V^k_\star$, it suffices to require that all the coefficients multiplying   $\|\X^{k+1}-\X^{k}\|^2$, $\|\X^{k}-\X^{k-1}\|^2$, $\|\Y^{k+1}-\Y^k\|^2$, $\|\Y^{k+1}-\Y^*\|^2$, and $\Delta  \mathcal{L}^k_\star$ in~\eqref{eq:Lyapunov dynamic} are non-negative. This is ensured by selecting, for a proper  $\zeta^k>0$,   $\alpha^k$ and   $\{\sigma^k\}_{k\geq 0}$  such that  the following hold for all $k\geq 0$  (note that~\eqref{eq:stepsize criterion 1}  simultaneously enforces non-negativity of the coefficients of $\|\X^{k+1}-\X^{k}\|^2$ and $\|\Y^{k+1}-\Y^k\|^2$; furthermore, since $\zeta^k>0$,  the coefficient of $\|\X^{k}-\X^{k-1}\|^2$, namely ${1}/{2}-\alpha^kL^k$, is also automatically non-negative):
\begin{subequations}
\label{eq:stepsize criteria}
    \begin{align}
        \label{eq:stepsize criterion 1}
\alpha^k\leq \min\left\{\frac{1}{2(L^k+\zeta^k)},\frac{\zeta^k}{\sigma^k}\right\};&\\
\label{eq:stepsize criterion 2}
(2+2\gamma^k)\alpha^k-2\gamma^{k+1}\alpha^{k+1}\geq& 0
;\\
\label{eq:stepsize criterion 3}
1/\sigma^k-1/\sigma^{k+1}\geq& 0.
\end{align}
\end{subequations}
Optimizing~\eqref{eq:stepsize criterion 1} over $\zeta^k>0$ (i.e., equating the two terms in the minimum), the largest admissible stepsize satisfying~\eqref{eq:stepsize criterion 1} is $\alpha^k= 1/(\sqrt{(L^k)^2+2\sigma^k}+L^k)$. Introducing  slack parameters $c_1,c_2\in(0,1]$ in~\eqref{eq:stepsize criterion 1} and~\eqref{eq:stepsize criterion 2},  respectively, and an additional    stepsize-control  policy 
$\pi^k:\mathbb{R}\rightarrow\mathbb{R}$ 
to add extra control on the stepsize selection policy,   conditions~\eqref{eq:stepsize criteria} are  satisfied by the following selection rule: 
\begin{equation}
 \label{eq:stepsize selection}\begin{aligned}\alpha^k=\min&\left\{\frac{1}{\sqrt{(L^k)^2+2\sigma^k/c_1}+L^k},\right.\\
 &\left.\sqrt{1+c_2\gamma^{k-1}}\alpha^{k-1},\pi^k(\alpha^{k-1})\right\}.\end{aligned}\end{equation}
 The stepsize-growth policy functions 
$\{\pi^k\}_{k\ge0}$   impose additional restrictions 
on possible increases of $\alpha^k$ when such control is needed 
for   convergence. When agents' losses are merely convex, this extra restriction is unnecessary and the term 
can be omitted, or one may simply take any admissible policy with 
$\pi^k(x)\ge x$ (see Assumption~\ref{ass:Pi local}). Our analysis shows that specific choices of $\pi^k$ (as in  
Assumption~\ref{ass:strong convx para}) are needed in the strongly convex setting to ensue   {\it linear} convergence with desirable rate scaling.   

The resulting line-search-free adaptive decentralized algorithm is summarized in Algorithm~\ref{alg:line-search-free}. Notice that the updates of the primal and dual variables (Step 3) require only one communication round/iteration within   neighboring nodes. 
\begin{algorithm}[htbp]
\centering
\resizebox{\columnwidth}{!}{%
  \begin{minipage}{\columnwidth}
\caption{\textbf{A}daptive \textbf{D}ecentralized \textbf{O}ptimization Algorithm with \textbf{L}ine-search-\textbf{F}ree Stepsize Selection (ADOLF)}
\label{alg:line-search-free}
  \noindent \textbf{Data:} $\X^{-1},\X^{0}\in\mathbb{R}^{m\times d}$,
   $\D^0=0$; $\alpha^0,\sigma^k>0$, $\forall k\geq 0$, $c_1,c_2\in(0,1]$, $\gamma^0=1$;   stepsize control policy   $\{\pi^k\}$;  $W=(1-c)I+c\widetilde{W}$, $c\in(0,1/2)$ and $\widetilde W\in \mathcal{W}_{\mathcal G}$.		
		\begin{algorithmic}[1]
        \State \texttt{(S.0)} Initialize: $$\begin{aligned}\D^1&=\D^0+\sigma^0\alpha^0(I-W)((1+\gamma^0)\X^0-\gamma^0\X^{-1}); \\\X^1&=\X^0-\alpha^0(\nabla F(\X^0)+\D^1);\end{aligned}$$
        \State \textbf{for} $k=1,2,\cdots \textbf{do}$
        \State \texttt{(S.1) Global scalar average: }Calculate$$L^k=\sqrt{\frac{\sum_{i\in[m]}\|\nabla f_i(x_i^k)-\nabla f_i(x_i^{k-1})\|^2}{\sum_{i\in[m]}\|x_i^k-x_i^{k-1}\|^2}};$$
           \State \texttt{(S.2) Stepsize selection:} Choose $\alpha^k$ according to~\eqref{eq:stepsize selection} and take $\gamma^k=\alpha^k/\alpha^{k-1}$;
			\State \texttt{(S.3) Dual and primal             update:}
         $$
             \begin{aligned}
\D^{k+1}&=\D^k+\sigma^k\alpha^k(I-W)((1+\gamma^{k})\X^k-\gamma^{k}\X^{k-1});\\
    \X^{k+1}&=\X^k-\alpha^k(\nabla F(\X^k)+\D^{k+1}). 
             \end{aligned}   
         $$
\end{algorithmic}
 \end{minipage}
    }
\end{algorithm}

\textit{On the global scalar average.}  Computing $L^k$ in (\texttt{S.1}) requires a global scalar average across the network.  In modern wireless mesh networks this can be realized over a low-rate, long-range interface such as LoRa~\cite{janssen2020lora} (broadcasting a single scalar per node in one hop), while vector transmissions  use a high-rate interface such as WiFi. Alternatively, the global scalar average can be replaced by a global-max consensus step at the cost of a slightly smaller stepsize. A fully local variant that removes any global communication requirement is presented in Section~\ref{sec:local-min}.
\vspace{-2mm}
\subsection{Convergence guarantees}
\label{sec:convergence}
 This subsection provides the convergence results for Algorithm~\ref{alg:line-search-free}. We postulate the following assumption.  
\begin{assumption}
    \label{ass:convex para}
   The sequence $\{\sigma^k\}_{k\geq 0}$  is monotonically nondecreasing and $0<\sigma^k\leq \overline{\sigma}$, for  some $0<\overline{\sigma}<\infty$ and all $k\geq 0$.    The   stepsize-growth policy functions    $\{\pi^k\}$, $\pi^k:\mathbb R\to \mathbb R$, satisfy   $\pi^k(x)\geq x$, for all   $x\geq 0$. 
\end{assumption}
\subsubsection{Convex $f_i$'s$-$Sublinear convergence rate}
\label{subsec:sublinear convergence}
 For any given saddle point  $(\X^*,\Y^*)$   of~\eqref{eq:lagarangian line-search-free}, let   
\begin{equation*}
    \mathcal{M}_*(\X):= F(\X)-F(\X^*)+\langle \PL\Y^*,\X-\X^*\rangle+\|\PL\X\|^2\geq 0.
\end{equation*}
Notice that $\mathcal{M}_*(\X)=0$ iff $\PL\X=0$ and $F(\X)=F(\X^*)$, meaning that $\X$ is a solution of Problem~\eqref{eq:primal}; furthermore,  $\mathcal{M}_*(\bullet)$ is continuous; hence, it is  a valid measure of optimality. Finally,  define the ergodic sequence: 
\begin{equation*}
      \overline{\X}^k:=\frac{1}{\theta^k}\sum_{t=0}^{k-1}\gamma^t\X^t,\quad\text{with}\quad \theta^k=\sum_{t=0}^{k-1}\gamma^t,\quad, \forall k\geq 1.
\end{equation*}
Sublinear convergence of Algorithm~~\ref{alg:line-search-free}  is summarized next. 
\begin{theorem}[sublinear rate]
\label{thm:convergence of iterates}
   Let $\{(\X^k,\D^k)\}_{k\geq 0}$ be the sequence  generated by Algorithm~\ref{alg:line-search-free} under Assumptions~\ref{ass:function}-\ref{ass:convex para}, with $\mu=0$ and  $c_1\in(0,1)$; and let $\{\Y^k\}_{k\geq 0}$ be  the auxiliary sequence     defined in~\eqref{eq:auxiliary update}.    Then, the following hold:
   
   \textbf{(a)} $\texttt{sequence convergence}$:  Sequences $\{\X^k\}$, $\{\D^k\}, $$\{\Y^k\}$ globally converge:   
        \begin{equation}
        \label{eq:convergence of iterates}
\X^k\rightarrow\X^{\star},\quad \Y^k\rightarrow\Y^{\star},\quad \D^k\rightarrow \PL\Y^{\star}, 
        \end{equation} where  $(\X^\star,\Y^\star)$ is  a saddle points of~\eqref{eq:lagarangian line-search-free};
        
        \textbf{(b)} $\texttt{sublinear rate}$:  Setting   $c_2\in(0,1)$,   it holds that 
        \begin{equation}
     \label{eq:sublinear rate}
\mathcal{M}_\star(\overline{\X}^k)\leq \varepsilon,\quad \text{for all }k\geq N_{\varepsilon}:=\mathcal{O}\left({\widetilde{L}}/{\varepsilon}\right),
 \end{equation} 
 where $\widetilde{L}$ is the  Lipschitz constant of $\nabla F$ restricted to the convex hull of the compact set $\cup_{t=0}^{\infty}[\X^{t-1},\X^t]$.
\end{theorem}
\begin{proof}
    See Appendix.~\ref{appx:proof of theorem 1}.
\end{proof}
 {Theorem~\ref{thm:convergence of iterates} certifies {\it sequence} convergence of   Algorithm~\ref{alg:line-search-free}   sublinear convergence rate of the merit function $\mathcal M_\star$, under  only \textit{local} smoothness of the $f_i$'s. This is a major departure from the  existing nonadaptive decentralized  methods (see, e.g., ~\cite{nedic2018network,sayed2014adaptation} and references therein), which are guaranteed to converge only under  {\it global} smoothness  of $f_i$'s. Furthermore,  they achieve a $\mathcal{O}(L/\varepsilon)$ complexity whereas  Algorithm~\ref{alg:line-search-free}'s rate  only depends on the {\it restricted} Lipschitz constant $\widetilde{L}$, which can be significantly smaller than $L$. This  certifies faster convergence,     under weaker assumptions.   Our rate improves also upon adaptive  decentralized  methods~\cite{Chen2026AdaptivePS,nazari2022dadam,chen2023convergence,li2024problem}, and matches the order of adaptive line-search-based  methods~\cite{chen2026adaptive,kuruzov2025adaptive,kuruzovachieving,chen2025parameter,Xu2025AnAP} while being line-search free. 

\subsubsection{Strongly convex $f_i$'s$-$Linear convergence}
\label{subsec:linear}
 We establish now linear convergence of   Algorithm~\ref{alg:line-search-free} when the agents' losses are locally strongly convex ($\mu>0$). 
 
 We find the following choice of $\{\sigma^k\}_{k\geq 0}$ and  $\{\pi^k(\bullet)\}$ convenient to achieve the desired rate scaling.    
 \begin{assumption}
     \label{ass:strong convx para}
      The sequence $\{\sigma^k\}_{k\geq 0}$  is  chosen as $\sigma^k=\sigma/(\alpha^k)^2$, for some $\sigma\in(0,c_1/2)$  and all $k\geq 0$. The functions $\pi^k:\mathbb R \to \mathbb R$ are selected as  $\pi^k(x)\leq ((k+\beta_1)/(k+1))^{\beta_2} x$, for some $\beta_1,\beta_2>0$ and any $k=0,1,\ldots .$
 \end{assumption}
 
 \begin{remark}
 It's easy to show that, under   Assumption~\ref{ass:strong convx para}, ~\eqref{eq:stepsize selection} can be equivalently rewritten as
 \begin{equation}
 \label{eq:stepsize selection linear implementation}
     \alpha^k=\min\left\{\frac{1/2-\sigma/c_1}{L^k},\sqrt{1+c_2\gamma^{k-1}}\alpha^{k-1},\pi^k(\alpha^{k-1})\right\},
 \end{equation}
 which can be computed using only local information. 
\end{remark}

\begin{theorem}[linear convergence]
    \label{thm:contraction}
 Let $\{(\X^k,\D^k)\}_{k\geq 0}$ be the sequence  generated by Algorithm~\ref{alg:line-search-free} under Assumptions~\ref{ass:function}-\ref{ass:W} and~\ref{ass:strong convx para}, with $\mu>0$ and  $c_1,c_2\in(0,1)$; and let $\{\Y^k\}_{k\geq 0}$ be  the auxiliary sequence     defined in~\eqref{eq:auxiliary update}.    Then,     

    \textbf{(a)} $\texttt{global sequence convergence}$: 
        \begin{equation}
        \label{eq:convergence of iterates linear}
\X^k\rightarrow\X^\star,\quad \Y^k\rightarrow\Y^\star,\quad \D^k\rightarrow \PL\Y^\star,
        \end{equation} where  $(\X^\star,\Y^\star)$ is   the saddle points $(\X^\star,\Y^\star)$  of~\eqref{eq:lagarangian line-search-free};

 \textbf{(b)} $\texttt{linear rate}$:    It holds that
 \begin{equation}
    \label{eq:R linear}
   \|\X^k-\X^*\|^2\leq \varepsilon,\quad  \forall k\geq N_{\varepsilon}:=\mathcal{O}\left(\frac{\widetilde{\kappa}}{1-\lambda_2(W)}\log (1/\varepsilon)\right),
\end{equation}
where $\widetilde{\kappa}:=\widetilde{L}/\widetilde{\mu}$, and    $\widetilde{\mu}$ is the  strong convexity parameter of $F$ restricted to  the {the convex hull of the compact set $\cup_{t=0}^{\infty}[\X^{t-1},\X^t]$.}\end{theorem}
\begin{proof}
See Appendix~\ref{appx:proof of theorem 2}.
\end{proof}

    As in the convex case, the linear rate  in~\eqref{eq:R linear} compares favorably with that of  nonadaptive methods; it  depends only on the {\it restricted} condition number $\widetilde{\kappa}$, generally smaller than the global one, $\kappa=L/\mu$. It also matches the rates of adaptive line-search methods~\cite{chen2026adaptive,kuruzov2025adaptive,kuruzovachieving,chen2025parameter}   while avoiding   line-search steps, thereby  reducing the computation  cost. 

 \section{Beyond global  scalar averages }
 \label{sec:local-min}
In this section, we introduce a fully local variant of
Algorithm~\ref{alg:line-search-free}, requiring  
neighbor-to-neighbor communications also for the scalar stepsize information.  The proofs for this section can be found in Appendix~\ref{appx:section 3}

   For each agent $i\in[m]$ and iteration $k\geq0$, let
$\alpha_i^k$, $\sigma_i^k$, and $\gamma_i^k$ denote the
local counterparts of $\alpha^k$, $\sigma^k$, and $\gamma^k$
in Algorithm~\ref{alg:line-search-free}. To simplify the presentation, we further denote by $\Lambda^k$, $\Sigma^k$ and $\Gamma^k$   the diagonal matrices of $\{\alpha^k_i\}_{i\in[m]}$, $\{\sigma^k_i\}_{i\in[m]}$ and $\{\gamma^k_i\}_{i\in[m]}$ respectively.   With   those local values, the primal dual variable are now written as
\begin{subequations}
\label{eq:primal dual update local}
             \begin{align}
             \label{eq:dual update_local}
\D^{k+1}&=\D^k+(I-W)\Sigma^k\Lambda^k((I+\Gamma^k)\X^k-\Gamma^{k}\X^{k-1});
\\
\label{eq:primal update_local}   
    \X^{k+1}&=\X^k-\Lambda^k(\nabla F(\X^k)+\D^{k+1}).  
             \end{align}   
            \end{subequations}
Supported by the convergence analysis in Sec.~\ref{sec:problem}, ideally we would like to  select the stepsize according to  \eqref{eq:stepsize selection}; which however is not compatible to neighboring communications. Next we put forth a  procedure that satisfies \eqref{eq:stepsize selection} almost all iterations while requiring only local communication.  

At iteration $k\geq1$, agent $i$ first computes the local curvature estimate (the local counterpart of \eqref{eq:L^k})
\[
L^k_i:=
\frac{\|\nabla f_i(x_i^k)-\nabla f_i(x_i^{k-1})\|}{\|x_i^k-x_i^{k-1}\|},
\]
with the convention $L^k_i=0$ if
$x_i^k=x_i^{k-1}$, and forms the curvature-based candidate
\begin{equation}
\label{eq:local curvature}
\hat{\alpha}_{i}^k
:=
\left(
\sqrt{(L^k_i)^2+2\sigma_i^k/c_1}+L^k_i
\right)^{-1}.
\end{equation}
Mimicking   \eqref{eq:stepsize selection}, one would   
set $\alpha_i^k=\min_{i\in[m]}\widetilde{\alpha}_i^k$ and $\gamma_i^k= {\alpha_i^k}/{\alpha_i^{k-1}}$, with  
\begin{equation}
    \label{eq:tilde alpha 1}
    \widetilde{\alpha}_i^k\leq\min\left\{\hat{\alpha}_i^k,\sqrt{1+c_2\gamma_i^{k-1}}\alpha_i^{k-1},\pi^k(\alpha_i^{k-1})\right\}.
\end{equation}
 However, this would call for   a global min-consensus step.  As proxy, we replace the global min therein with a local one and adjust the definition of $\gamma^k$ as follows:  
at each $k\geq 0$, given $\widetilde{\alpha}_i^k$ defined in~\eqref{eq:tilde alpha 1}, let
 \begin{equation}
     \label{eq:stepsize selection local 1}
\alpha_i^k=\min_{j\in\mathcal{N}_i}\widetilde{\alpha}_j^k\quad \text{and}\quad \gamma_i^k=\frac{\widetilde{\alpha}_i^k}{\alpha_i^{k-1}},\quad i\in [m].
 \end{equation}
 This choice is implementable by neighboring communications. Furthermore,   under  suitable  $\{\pi^k\}$ and $\sigma_i^k$'s,  as stated in  Assumption~\ref{ass:Pi local} below,  (\ref{eq:stepsize selection local 1}) achieves {\it global} consensus in {\it finite} time, still preserving adaptivity in the consensual stepsize value.  Lemma~\ref{lemma:conditioned stepsize consistency} below captures this property. 
 \begin{assumption}
\label{ass:Pi local} 
      For each $i\in[m]$, the sequence $\{\sigma^{k,i}\}_{k\geq 0}$  is monotonically nondecreasing and $0<\sigma^{k,i}\leq \overline{\sigma}$, for  some $0<\overline{\sigma}<\infty$ and all $i\in[m]$ and $k\geq 0$. The stepsize policy correction sequence $\{\pi^k\}$, $\pi^k:\mathbb{R}\rightarrow\mathbb{R}$, satisfies
      \begin{equation}
      \label{eq:pi local}
          \pi^k(x)\leq x  +\delta^k,\quad, \forall k\geq 0, 
      \end{equation} 
  where  $\{\delta^k\}$ is any   summable, nonnegative sequence.
\end{assumption}
 \begin{lemma}
    \label{lemma:conditioned stepsize consistency}
    Let $\{(\X^k,\D^k)\}$ be the sequence of iterates generated by~\eqref{eq:primal dual update local} under Assumption~\ref{ass:Pi local}, with $\{\Lambda^k\}$, $\{\Gamma^k\}$ selected by~\eqref{eq:stepsize selection local 1}. 
    If the set $\mathcal{K}:=\{k:\exists i\in[m], \hat{\alpha}_i^k\leq \pi^k(\alpha_i^{k-1})\}$ is finite and $\{\X^k\}$ is bounded, then there exists $K>0$  such that, for any $k\geq K$, it holds
      $$\Lambda^k=\alpha^k I \quad\text{and}\quad \Gamma^k=\gamma^k I,\quad\text{with}\quad  \gamma^k=\frac{\alpha^k}{\alpha^{k-1}},$$
       and    $\alpha^k$ satifying~\eqref{eq:stepsize selection}, with $\sigma^k:=\min_{i\in[m]}\sigma_i^k$. 
\end{lemma}
It remains to ensure $\hat{\alpha}_i^k\leq \pi^k(\alpha_i^{k-1})$ occurs finitely often. To this end, with $\eta\in(0,1)$, we adopt the strategy to enforce a sufficient decrease of $\alpha_i^k$, when $\hat{\alpha}_i^k\leq \pi^k(\alpha_i^{k-1})$ by taking
 \begin{equation}
\label{eq:tilde alpha}
\widetilde{\alpha}_{i}^k :=
\begin{cases}
\min\{\eta\alpha^{k-1}_i, \hat{\alpha}^k_{i}\},
  \quad\text{if }\hat{\alpha}^k_{i} \leq \pi^k(\alpha^{k-1}_i); \\
\min\Big\{\pi^k(\alpha^{k-1}_i),
 \sqrt{1+c_2\gamma^{k-1}_i}\,\alpha^{k-1}_i
\Big\},
\quad \text{else}.
\end{cases}
\end{equation}
In fact, $\widetilde{\alpha}_i^k$ in~\eqref{eq:tilde alpha} satisfies~\eqref{eq:tilde alpha 1} and certifies Lemma~\ref{lemma:finite K}.
\begin{lemma}
    \label{lemma:finite K}
    Instate   Lemma~\ref{lemma:conditioned stepsize consistency}, but with  $\{\widetilde{\alpha}_i^k\}$ chosen according to~\eqref{eq:tilde alpha}. Then $\{\alpha_i^k\}$ is bounded from below and $|\mathcal{K}|<\infty$.
\end{lemma}
The algorithm, based on the above stepsize selection,  is summarized   in  Algorithm~\ref{alg:line-search-free_local}. 
 \begin{algorithm}[ht]
\centering
\resizebox{.9\columnwidth}{!}{%
  \begin{minipage}{\columnwidth}
    \caption{ADOLF-local}
\label{alg:line-search-free_local}
  \noindent \textbf{Data:} $\X^{-1},\X^{0}\in\mathbb{R}^{m\times d}$,
   $\D^0=0$; $\Lambda^0, \Sigma^k\succ0, \forall k\geq 0, \Gamma^0=I;$ $c_1,c_2\in(0,1]$ $\eta\in(0,1)$; stepsize control policy $\{\pi^k\}$;    $W=(1-c)I+c\widetilde{W}$, $c\in(0,1/2)$.		
		\begin{algorithmic}[1]
        \State \texttt{(S.0)} Initialize:
        $$\begin{aligned}
            \D^1&=\D^0+(I-W)\Sigma^0\Lambda^0((I+\Gamma^0)\X^0-\Gamma^0\X^{-1});\\
            \X^1&=\X^0+\Lambda^0(\nabla F(\X^0)+\D^1);
        \end{aligned}$$
        \State \textbf{for} $k=1,2,\cdots \textbf{do}$
        \State \texttt{(S.2) Local stepsize selection: }Select $\widetilde{\alpha}_i^k$ according to~\eqref{eq:tilde alpha}.
           \State \texttt{(S.3) Local min-consensus: }Select $\alpha_i^k$ and $\gamma_i^k$ according to~\eqref{eq:stepsize selection local 1}.

			\State \texttt{(S.4) Dual and Primal             update: }Update $\D^{k+1}$ and $\X^{k+1}$ according to~\eqref{eq:primal dual update local}.
\end{algorithmic}
 \end{minipage}
    }
\end{algorithm}

 Theorem~\ref{thm:iterates convergence_local} summarizes sublinear convergence of Algorithm~\ref{alg:line-search-free_local} when  the $f_i$'s are merely convex, where we defined the following  ergodic sequences: given $K\geq 0$, let 
 {\begin{equation}
    \label{eq:ergodic sequence local}
    \!\!\!\overline{\X}_K^k:=\frac{1}{\theta_K^k}\sum_{t=K}^{k-1}\gamma^t\X^t,\quad \theta_K^k=\sum_{t=K}^{k-1}\gamma^t,\quad  k\geq K+1.
\end{equation}
}
\begin{theorem}[sublinear rate]
\label{thm:iterates convergence_local}
      Let $\{(\X^k,\D^k)\}$ be the sequence of iterates generated by Algorithm~\ref{alg:line-search-free_local}  under Assumption~\ref{ass:function}-\ref{ass:W}, and~\ref{ass:Pi local}, with $\mu=0$ and $c_1\in(0,1)$; further assume that  
      $\{\X^k\}$ is bounded. Then, the following hold:
   
   \textbf{(a)} \textit{global sequence convergence: }  
        \begin{equation}
        \label{eq:convergence of iterates local}
\X^k\rightarrow\X^{\infty},\quad  \D^k\rightarrow \PL\Y^{\infty}:=\D^{\infty},
        \end{equation} where  $(\X^{\infty},\Y^{\infty})$ is a saddle points $(\X^\star,\Y^\star)$  of~\eqref{eq:lagarangian line-search-free};
        
        \textbf{(b)} \textit{sublinear rate: } Setting $c_2\in(0,1)$,  $\exists K>0$ such that  
 \begin{equation}
     \label{eq:sublinear rate local}
 {\mathcal{M}_*(\overline{\X}_K^k)}\leq \varepsilon,\quad \text{for all }k\geq N_{\varepsilon}:=\mathcal{O}\left(K+\widetilde{L}_K/\varepsilon\right),
 \end{equation} 
 where $\widetilde{L}_K$ is the Lipschitz constant of $\nabla F$ over the convex hull of compact set $\cup_{t=K}^{\infty}[\X^{t-1},\X^t]$.
\end{theorem}
Similar to Algorithm~\ref{alg:line-search-free}, Algorithm~\ref{alg:line-search-free_local} achieves linear convergence with a favorable rate under another choice of $\{\Sigma^k\}$ when $\mu>0$ ($f_i$'s strongly convex). 
 \begin{assumption}
     \label{ass:strong convx Sigma local}
      The sequence $\{\Sigma^k\}_{k\geq 0}$  is selected such that $\sigma^k_i=\sigma/(\alpha^k_i)^2$ for some $\sigma\in(0,c_1/2)$ for any $k\geq 0$ and $\{\pi^k\}$ is selected according to~\eqref{eq:pi local}.
 \end{assumption}
\begin{theorem}[linear convergence]
    \label{thm:contraction local}
Let $\{(\X^k,\D^k)\}_{k\geq 0}$ be the sequence  generated by Algorithm~\ref{alg:line-search-free_local} under Assumptions~\ref{ass:function}-\ref{ass:W}, and~\ref{ass:strong convx Sigma local}, with $\mu>0$ and and  $c_1,c_2\in(0,1)$. 
 Let $(\X^*,\Y^*)$ be the unique saddle point of~\eqref{eq:lagarangian line-search-free}, then if the sequence of iterates $\{\X^k\}$ is bounded, the following hold:

  {\textbf{(a)} \textit{Global sequence convergence: }
        \begin{equation}
        \label{eq:convergence of iterates linear local}
\X^k\rightarrow\X^*,\quad \D^k\rightarrow \PL\Y^*:=\D^*,
        \end{equation} }
\textbf{(b)} \textit{Linear rate: } There exists a finite   $K>0$ such that $\|\X^k-\X^*\|^2\leq \varepsilon$, for all $k\geq N_\varepsilon$, with
\begin{equation}
    \label{eq:R linear local}
N_{\varepsilon}=\mathcal{O}\left(K+\frac{\widetilde{\kappa}_{K}}{1-\lambda_2(W)}\log(1/\varepsilon)\right),
\end{equation}
where $\widetilde{\kappa}_{K}:=\widetilde{L}_{K}/\widetilde{\mu}_{K}$ and $\widetilde{\mu}_K$ is the  strong convexity parameter of $F$ restricted to  the {the convex hull of the compact set $\cup_{t=K}^{\infty}[\X^{t-1},\X^t]$.}
\end{theorem}
 A comparison between Algorithms~\ref{alg:line-search-free} and
\ref{alg:line-search-free_local} shows the price paid for removing the
global scalar average:   In the local variant, agents may initially use
different stepsizes, inducing a {\it nonmonotone} decrease of the merit function   from the first
iteration. This mismatch is only transient, that is,   there exists a finite index $K$ such that, for all
$k\ge K$, the local method uses a common scalar stepsize (still changing), i.e.,
$\Lambda^k=\alpha^k I$ and $\Gamma^k=\gamma^k I$. The offset $K$ in Theorems~\ref{thm:iterates convergence_local} and
\ref{thm:contraction local} accounts exactly for this initial
nonhomogeneous phase.  After iteration $K$, the monotonically decrease of the merit function is restored. 

The boundedness assumption in
Theorems~\ref{thm:iterates convergence_local} and
\ref{thm:contraction local} is due to the fact that  one only assumed  local
smoothness, rather than the global one, commonly used
in decentralized optimization. Under local smoothness,  curvature
constants are uniformly bounded (hence the stepsizes)   only on bounded regions;  boundedness of
the generated trajectory ensures then that the restricted constants $\widetilde{L}$ and $\widetilde{\kappa}$
are uniformly bounded. This type of assumption is
standard in the  analyses under   local Lipschitz gradients. It is also
not a tuning requirement: the algorithm does not need to know any bound
on the iterates. If the standard global-smoothness assumption is postulated
instead, as in most existing decentralized methods, our analysis apply readily just  the global smoothness constants, and {\it the boundedness of the iterates
guaranteed}; hence Theorems~\ref{thm:iterates convergence_local} and
\ref{thm:contraction local} will be unconditional.  

\section{Numerical Results}
\label{sec:numeric result}
This section presents some preliminary numerical results, comparing Algorithms~\ref{alg:line-search-free} (ADOLF) and~\ref{alg:line-search-free_local} (ADOLF-local) with the following   benchmarks:   EXTRA~\cite{shi2015extra} (example of nonadaptive algorithm) and the recent adaptive proposals    DATOS and DATOS-local~\cite{chen2026adaptive}. EXTRA requires full knowledge of network and optimization parameters; we select the stepsize based on  a grid-search to identify the choice  ensuring fastest practical convergence. DATOS and DATOS-local perform a local line-search to select the stepsize, followed by a global and local-min consensus (communication) step, respectively. In ADOLF-local, we chose $\{\pi^k\}$ such that $\pi^k(x)=x+\frac{6}{\pi^2}\cdot\frac{1}{k^2}$, and $\eta=0.9$. No  $\{\pi^k\}$ is used  in ADOLF  for the case-study in Sec~\ref{subsec:logistic} whereas we chose  $\pi^k(x)=[(k+10)/(k+1)]x$ for ADOLF applied to the ridge regression problem in Sec~\ref{subsec:ridge}.

We simulate   line graphs  and Erdos-Renyi graphs with  edge-probability   $p=0.1$ and  $p=0.9$; there are    $m=20$ agents. The gossip weights  are the Metropolis-Hasting    \cite{Nedic_Olshevsky_Rabbat2018}. 

\subsection{Logistic regression}
\label{subsec:logistic}
Consider  the decentralized logistic regression problem, which is an instance of \eqref{eq:problem}, with 
\begin{equation*}\label{eq:logistic-scvx}f_i(x)=\frac{1}{n}\sum_{j=1}^{n}\log(1+\exp(-b_{ij}\cdot\langle x,a_{ij}\rangle)),\end{equation*}
where  $a_{ij}\in \mathbb{R}^{d}$ and $b_{ij}\in\{-1,1\}$. Data $\{(a_{ij},b_{ij})\}_{j=1}^n$ is  owned   by agent $i$. We use the MNIST dataset. 
The feature dimension is $d=784$. 

Figure~\ref{fig:general cvx}   plots the optimality gap  $({1}/{m})\sum_{i=1}^m f(x_i)-f^*$ versus the number of communications, achieved by all the algorihtms,  where $u$ is the objective function in~\eqref{eq:problem} and $u^*$ is its minimal value. 
The figures show  that both  proposed methods consistently outperform the EXTRA algorithm with fixed stepsize, with the advantage that do not require any user's intervention for the tuning of the stepsize. Moreover, our methods perform comparably with DATOS an DATOS-local, with the advantage that no extra line-search process. 

\subsection{Ridge regression}
\label{subsec:ridge}
We also report experiments solving a ridge regression problem,  a strongly convex instance of~\eqref{eq:problem},  with 
\begin{equation}
    f_i(x)=\frac{1}{n}\|A_ix-b_i\|^2+\frac{\gamma^i}{2}\|x\|^2,
\end{equation}
 where $(A_i,b_i)\in\mathbb{R}^{n\times d}\times \mathbb{R}^n$ are the data  owned  by agent $i$. The elements of $A_i$, $b_i$ are independently sampled from the
standard normal distribution. Here,   $n=20$ and $d=500$. We set $\gamma^i=0.1+(i-1)\times 0.1$, so that the smoothness parameter $L_i$ are different among each local function. 

 Figure~\ref{fig:general scvx} plots the optimality gap measured by $\|\X^k-\X^*\|^2$ versus the number of communications. Similar to the convex scenario, both proposed algorithms clearly outperform existing decentralized benchmarks. 

\begin{figure*}[htbp] 
    \centering
    \includegraphics[width=0.25\linewidth]{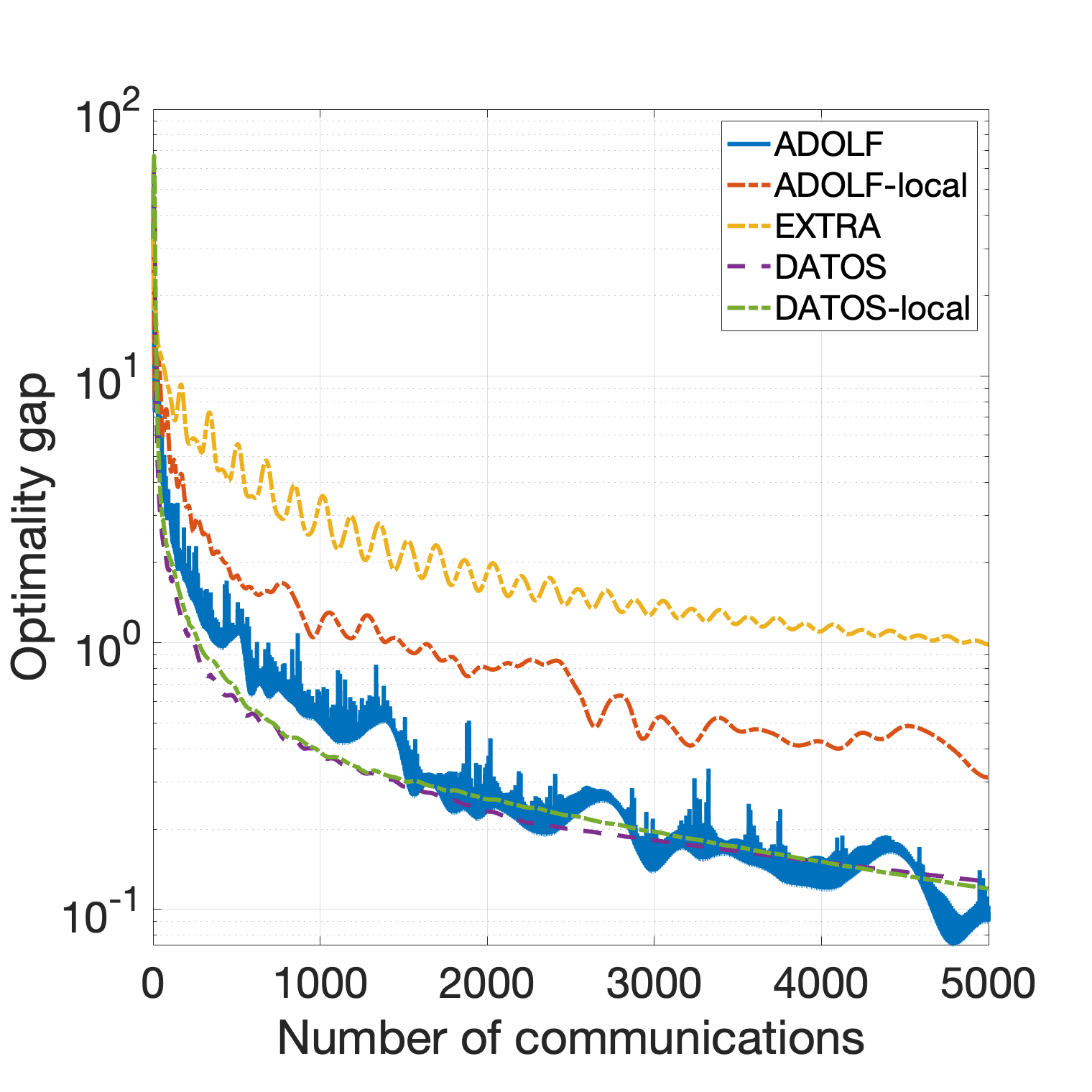}
    \hfill
    \includegraphics[width=0.25\linewidth]{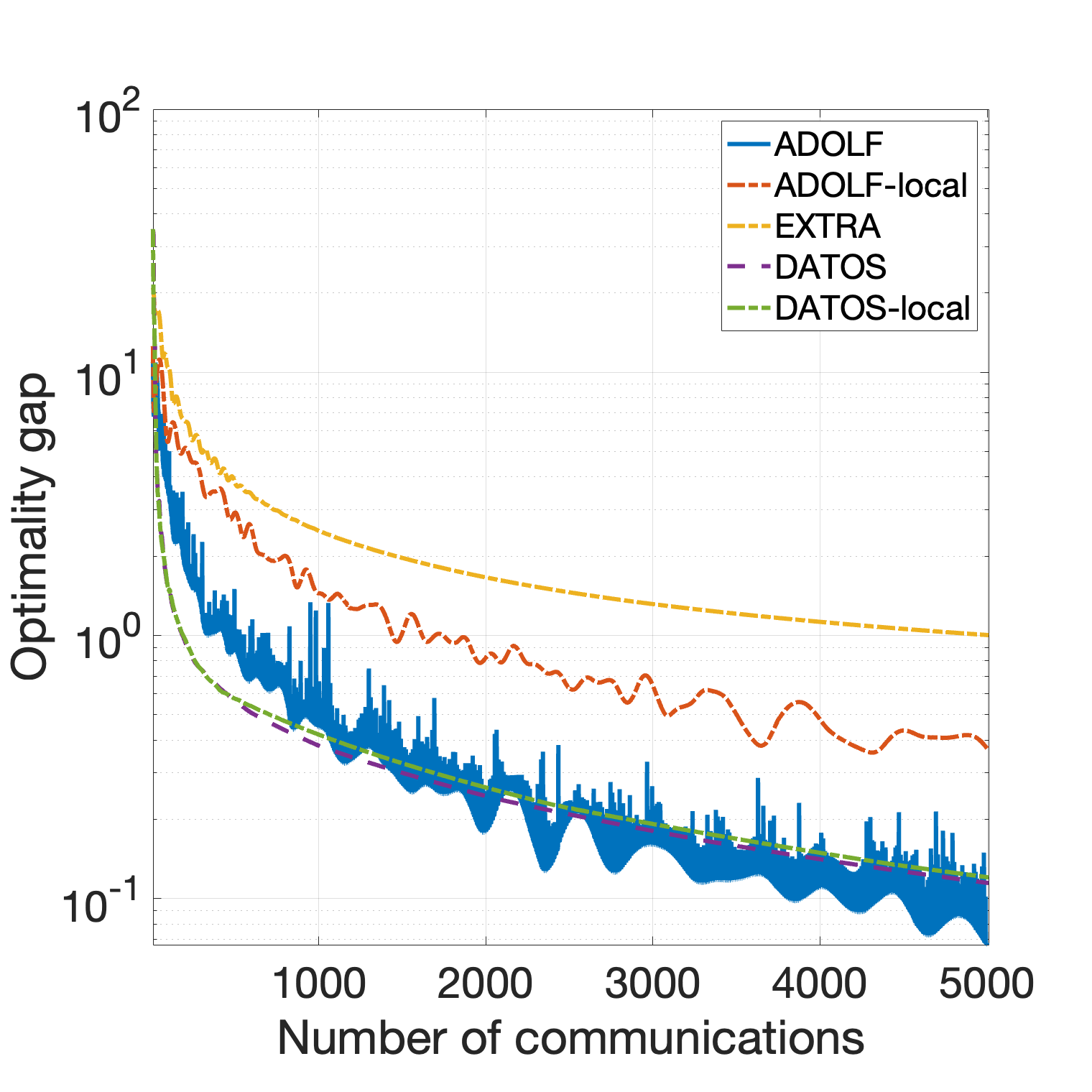}
    \hfill
    \includegraphics[width=0.25\linewidth]{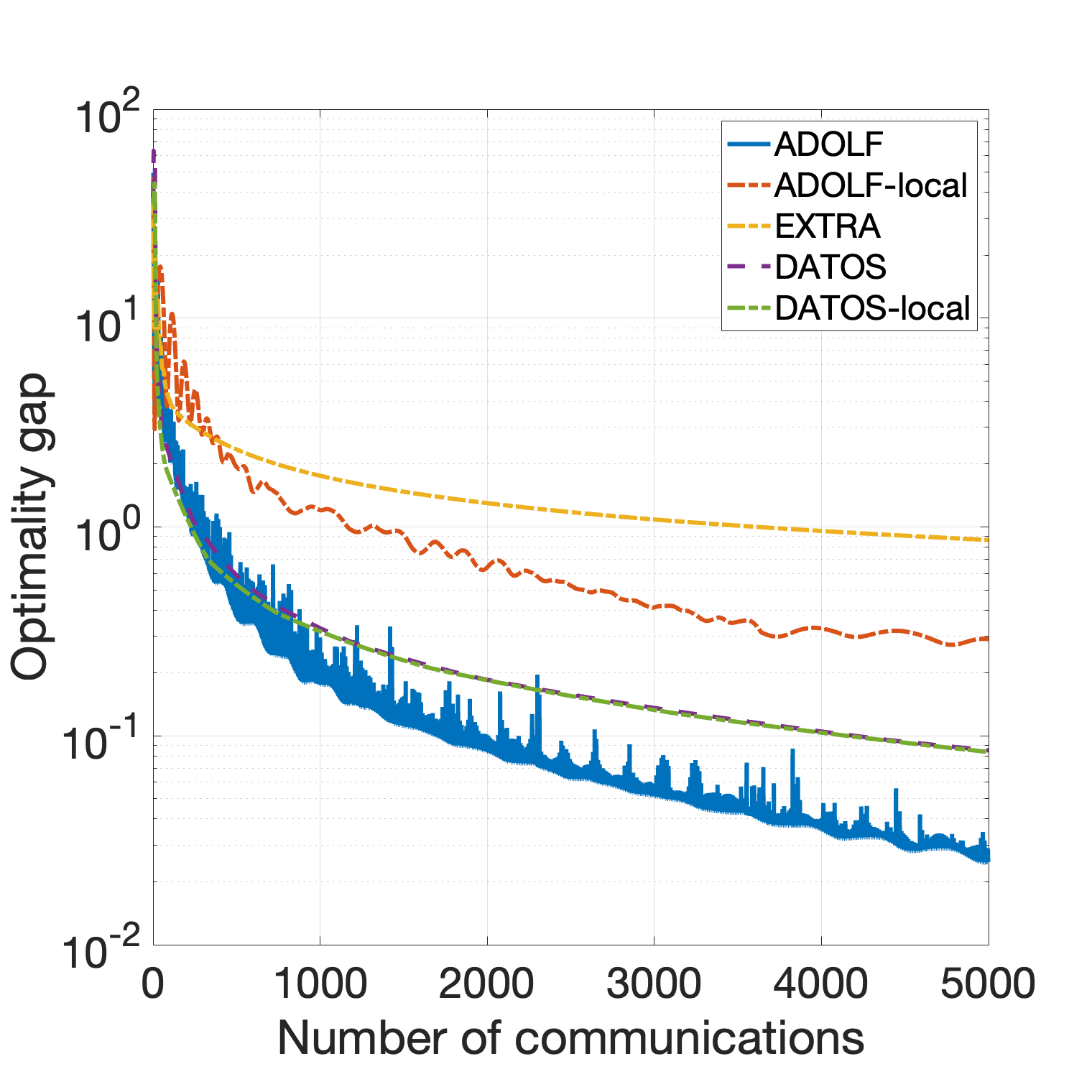}
    \caption{Logistic regression: ${\frac{1}{m}\sum_{i=1}^m f(x_i)-f(x^*)}$ v.s. \# communications.  Comparison of EXTRA, DATOS, DATOS-local, ADOLF and ADOLF-local on a line graph (left) and Erdos-Renyi graphs with different edge-probability:  $p=0.1$  (middle); and  $p=0.9$ (right).}\vspace{-0.45cm}

\label{fig:general cvx}
\end{figure*}
  \begin{figure*}[htbp]
    \centering
    \includegraphics[width=0.25\linewidth]{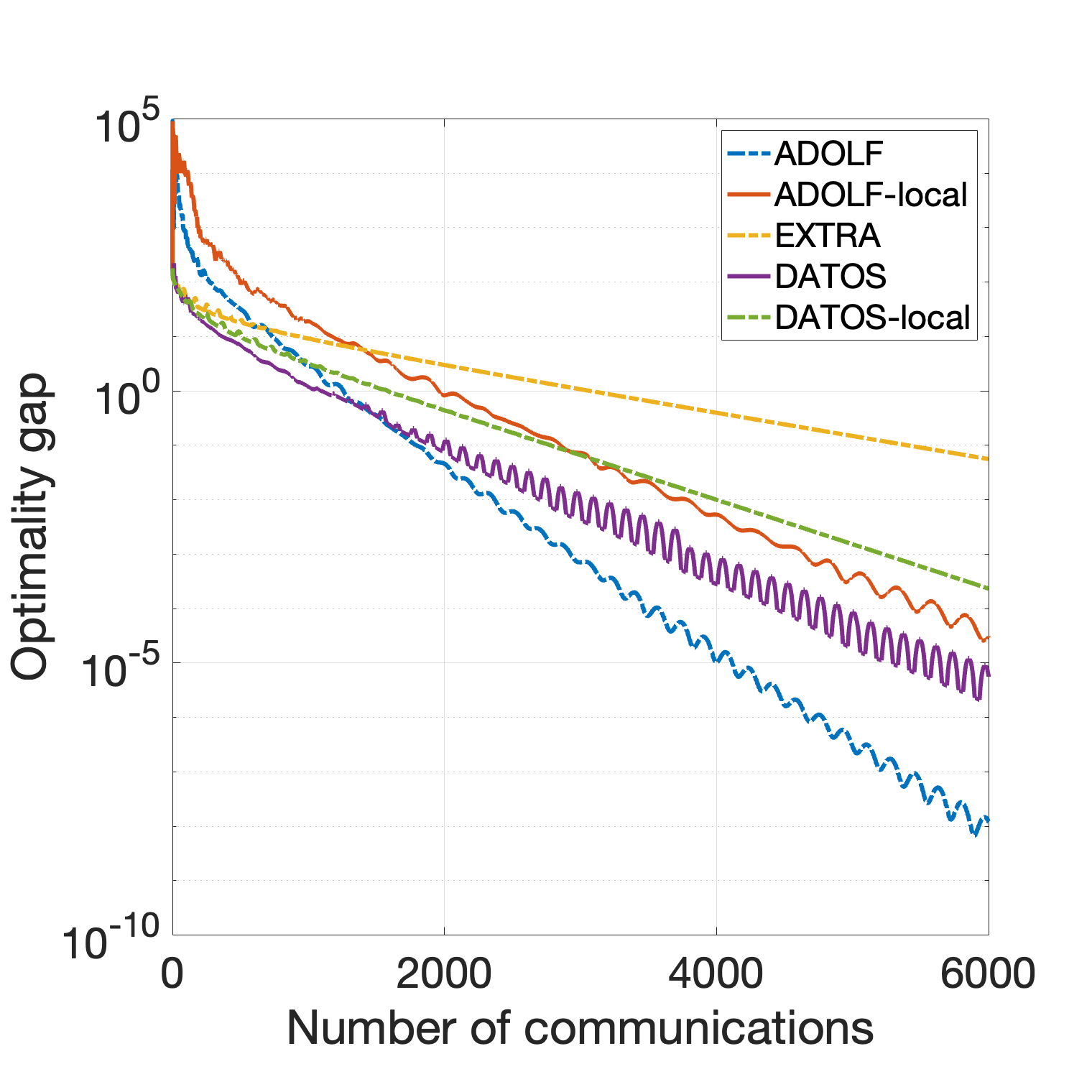}
    \hfill
    \includegraphics[width=0.25\linewidth]{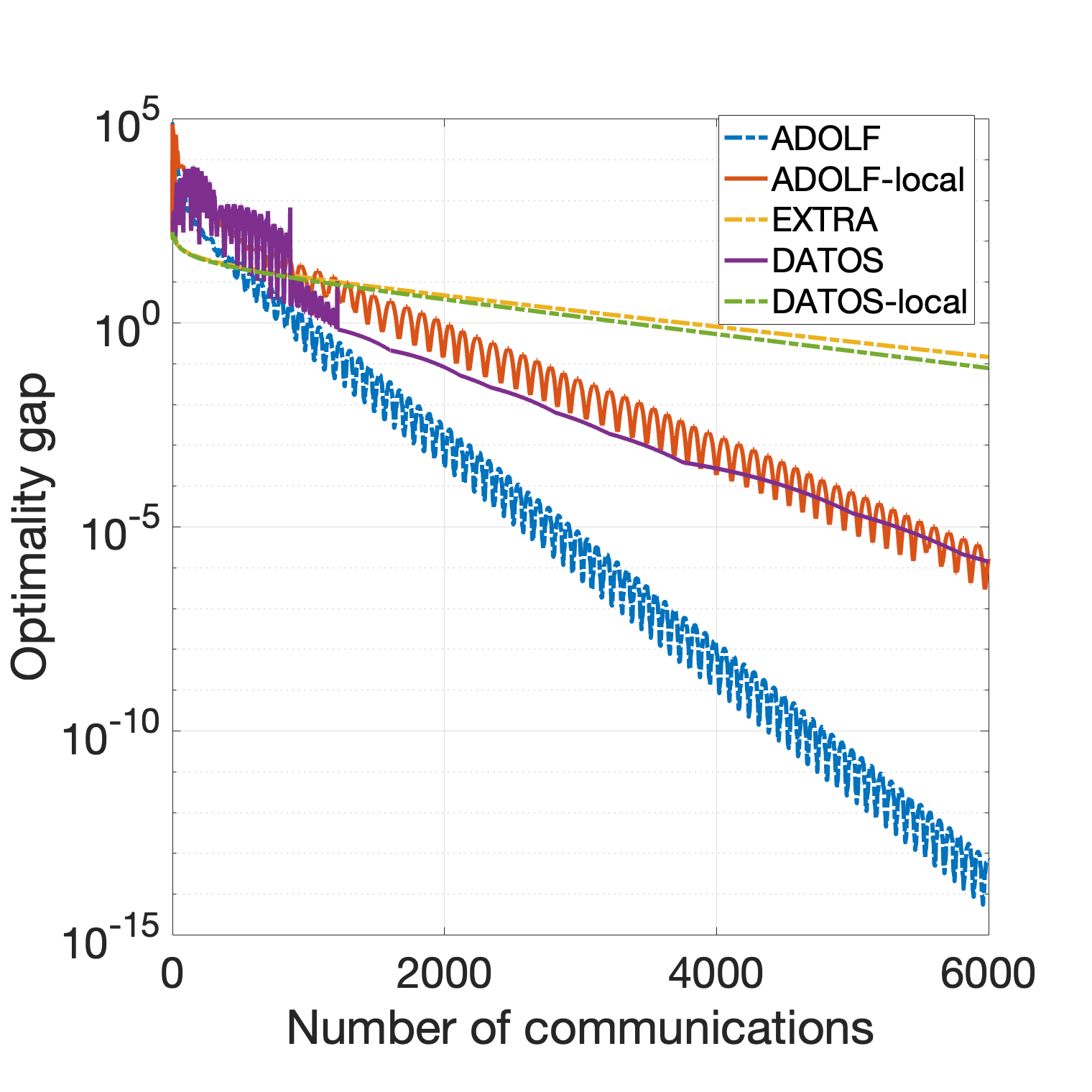}
    \hfill
    \includegraphics[width=0.25\linewidth]{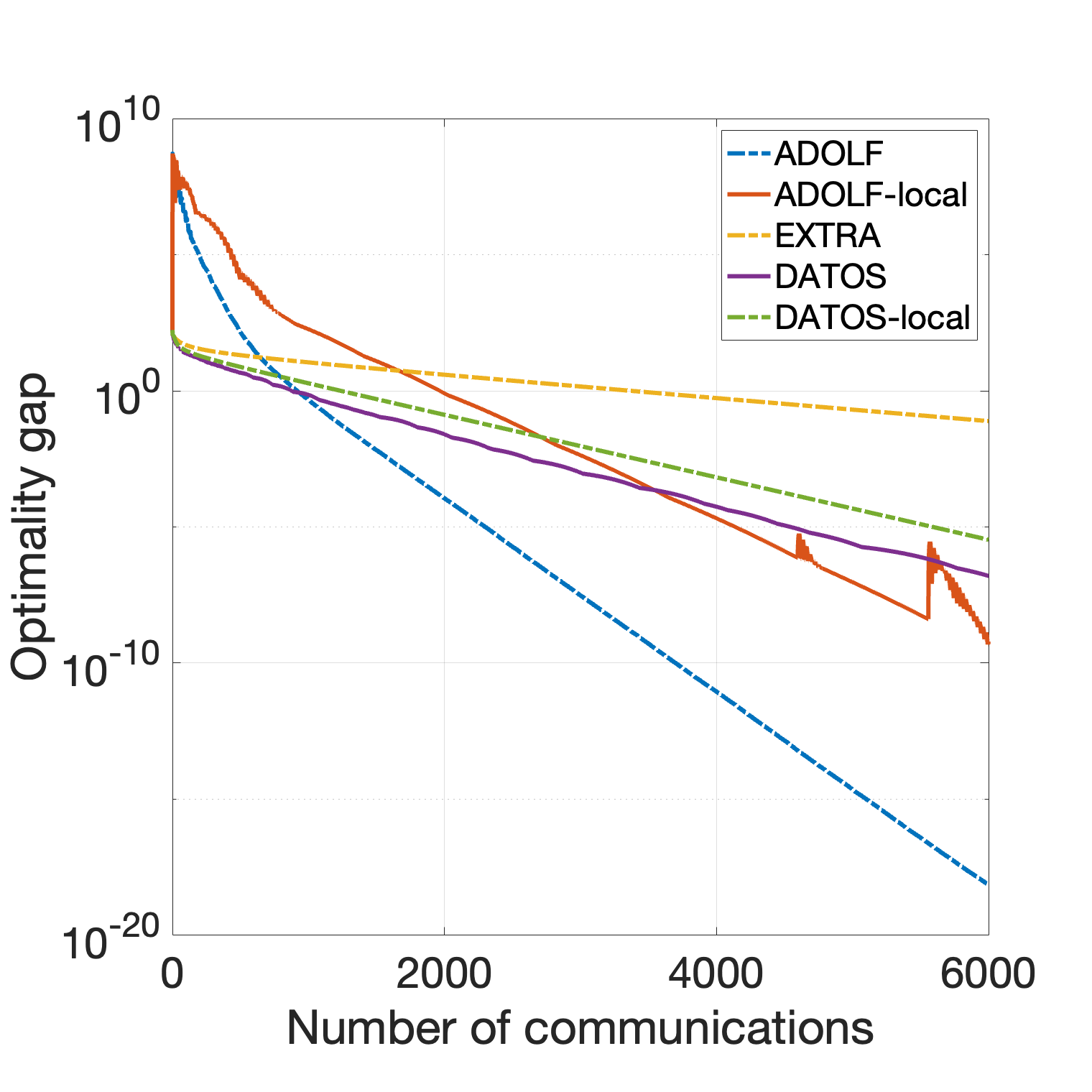}
\caption{Linear regression with $\ell_2$ regularization: ${\|\X^k-\X^*\|^2}$  v.s. \# communications. Comparison of EXTRA, DATOS, DATOS-local, ADOLF, and ADOLF-local on a line graph (left) and Erdos-Renyi graphs with different edge-probability:  $p=0.1$  (middle); and  $p=0.9$ (right).}\vspace{-0.45cm} 
\label{fig:general scvx}
\end{figure*}
\bibliographystyle{ieeetr-diss}
\bibliography{Bib}

\appendices
\section{Proof of Lemma~\ref{lemma:Lyapunov dynamic}}
\label{appx:proof of lemma lyapunov}
\begin{proof}
    Consider for $k\geq 0$, \begin{equation}
\begin{aligned}
\label{eq:primal first step}
    &\|\X^{k+1}-\X^*\|^2\\
    =&\|\X^k-\X^*\|^2-\|\X^{k+1}-\X^{k}\|^2\\
    &-2\alpha^k\langle \nabla F(\X^k)+\D^{k+1},\X^k-\X^*\rangle\\
    &+\underbrace{2\alpha^k\langle \nabla F(\X^k)+\D^{k+1},\X^k-\X^{k+1}\rangle}_{\texttt{term I}}.
\end{aligned}
    \end{equation}
    And further decompose $\texttt{term I}$, we have for $k\geq 0$,
    \begin{equation}
    \label{eq:term I}
        \begin{aligned}
            &\texttt{term I}\\
\overset{\eqref{eq:primal update}}{=}&2\alpha^k\langle\nabla F(\X^k)-\nabla F(\X^{k-1}),\X^k-\X^{k+1}\rangle\\
&+2\gamma^{k}\alpha^k\langle \X^{k-1}-\X^k,\nabla F(\X^k)+\D^{k+1}\rangle\\
&+2\alpha^k\langle \Y^{k+1}-\Y^k,\PL(\X^{k}-\X^{k+1})\rangle\\
\leq& \alpha^kL^k \left(\|\X^k-\X^{k-1}\|^2+\|\X^{k+1}-\X^{k}\|^2\right)\\
&+2\gamma^{k}\alpha^k\langle \X^{k-1}-\X^k,\nabla F(\X^k)+\D^{k+1}\rangle\\
+&\alpha^k\left(\zeta^k\|\PL(\X^{k+1}-\X^{k})\|^2+\frac{1}{\zeta^k}\|\Y^{k+1}-\Y^k\|^2\right),
        \end{aligned}
    \end{equation}
    where $\zeta^k>0$ and the last inequality comes from Young's inequality. Moreover, for $k\geq 0$, 
    \begin{equation}
        \label{eq:separation equality}
\begin{aligned}
            \|\PL(\X^{k+1}-\X^{k})\|^2=&\langle\X^{k+1}-\X^{k},\PL^2(\X^{k+1}-\X^{k})\rangle\\
            =&\|\X^{k+1}-\X^{k}\|^2-\|\X^{k+1}-\X^{k}\|^2_W
\end{aligned}
    \end{equation}
Replace~\eqref{eq:separation equality} into $\texttt{term I}$ and replace $\texttt{term I}$ into~\eqref{eq:primal first step}, we have for $k\geq 0$ that\begin{equation}
    \label{eq:primal last step}
    \begin{aligned}
        &\|\X^{k+1}-\X^*\|^2\\
        \leq&\|\X^k-\X^*\|^2-(1-\alpha^kL^k-\alpha^k\zeta^k)\|\X^{k+1}-\X^{k}\|^2\\
        &+\alpha^kL^k\|\X^{k}-\X^{k-1}\|^2+2\alpha^k\langle\Y^{k+1},\PL(\X^*-\X^k)\rangle\\
        &-\alpha^k\zeta^k\|\X^{k+1}-\X^{k}\|_W^2+2\gamma^{k}\alpha^k(F(\X^{k-1})-F(\X^*))\\
        &+2\gamma^{k}\alpha^k\langle \PL(\X^{k-1}-\X^k),\Y^{k+1}\rangle+\frac{\alpha^k}{\zeta^k}\|\Y^{k+1}-\Y^k\|^2\\
        &+(2+2\gamma^{k})\alpha^k(F(\X^*)-F(\X^k)),
    \end{aligned}
\end{equation}
Then, consider update of auxiliary variable, for $k\geq 0$,
\begin{equation}
    \label{eq:auxiliary first step}
    \begin{aligned}
        \|\Y^{k+1}-\Y^*\|^2=&\|\Y^k-\Y^*\|^2-\|\Y^k-\Y^{k+1}\|^2\\
        &+\underbrace{2\langle \Y^{k+1}-\Y^*,\Y^{k+1}-\Y^k\rangle}_{\texttt{term II}}.
    \end{aligned}
\end{equation}
By the update of auxiliary update~\eqref{eq:auxiliary update}, we have for $k\geq 0$,
\begin{equation}
    \label{eq:term II}
    \begin{aligned}
        \texttt{term II}=&
        2\sigma^k\alpha^k\langle\PL((1+\gamma^{k})\X^k-\gamma^{k}\X^{k-1}),\Y^{k+1}-\Y^{*}\rangle\\
        =& 
        2\sigma^k\gamma^{k}\alpha^k\langle\PL(\X^k-\X^{k-1}),\Y^{k+1}-\Y^{*}\rangle\\
        &+2\sigma^k\alpha^k\langle\PL\X^k,\Y^{k+1}-\Y^{*}\rangle.
    \end{aligned}
\end{equation}
Replace $\texttt{term II}$ into~\eqref{eq:auxiliary first step}, we have for $k\geq 0$ that,
\begin{equation}
    \label{eq:auxiliary last step}
    \begin{aligned}
   \|\Y^{k+1}-\Y^*\|^2=&\|\Y^k-\Y^*\|^2-\|\Y^k-\Y^{k+1}\|^2\\
   &
   +2\sigma^k\gamma^{k}\alpha^k\langle\PL(\X^k-\X^{k-1}),\Y^{k+1}-\Y^{*}\rangle\\
   &+2\sigma^k\alpha^k\langle\PL\X^k,\Y^{k+1}-\Y^{*}\rangle.
    \end{aligned}
\end{equation}
Take$~\eqref{eq:primal last step}+\frac{1}{\sigma^k}~\eqref{eq:auxiliary last step}$, it holds that for $k\geq 0$,
\begin{equation}
    \label{eq:combined first step}
    \begin{aligned}
        &\|\X^{k+1}-\X^*\|^2+\frac{1}{\sigma^k}\|\Y^{k+1}-\Y^*\|^2\\
        \leq& \|\X^{k}-\X^*\|^2+\frac{1}{\sigma^k}\|\Y^{k}-\Y^*\|^2\\
        &-(1-\alpha^kL^k-\alpha^k\zeta^k)\|\X^{k+1}-\X^{k}\|^2\\
        &+\alpha^kL^k\|\X^{k}-\X^{k-1}\|^2-\left(\frac{1}{\sigma^k}-\frac{\alpha^k}{\zeta^k}\right)\|\Y^{k+1}-\Y^k\|^2\\
        &-\alpha^k\zeta^k\|\X^{k+1}-\X^{k}\|_W^2-(2+2\gamma^{k})\alpha^k\langle\X^k-\X^*,\PL\Y^*\rangle\\
        &-2\gamma^{k}\alpha^k\langle\X^*-\X^{k-1},\PL\Y^*\rangle+2\alpha^k\langle \Y^{k+1}-\Y^*,\PL\X^*\rangle\\
        &
        +2\gamma^{k}\alpha^k(F(\X^{k-1})-F(\X^*))\\
        &+(2+2\gamma^{k})\alpha^k(F(\X^*)-F(\X^k)).
    \end{aligned}
\end{equation}
Replace the definitions~\eqref{eq:L_p gap} into~\eqref{eq:combined first step} and by $\PL\X^*=0$,
~\eqref{eq:Lyapunov dynamic} can be obtained by adding the term $\frac{1}{2}\|\X^{k+1}-\X^k\|^2+\frac{1}{2}\|\X^k-\X^{k-1}\|^2+2\gamma^{k+1}\alpha^{k+1}\Delta \mathcal{L}_*^k+(1/\sigma^{k+1}-1/\sigma^k)\|\Y^{k+1}-\Y^k\|^2$ on both sides of the resulting inequality and reordering terms.
\end{proof}
\section{Proof of theorem~\ref{thm:convergence of iterates}}
\label{appx:proof of theorem 1}
We display several useful lemma to prove Theorem~\ref{thm:convergence of iterates}. Let's start from the dynamic of Lyapunov function. Based on the stepsize selection policy~\eqref{eq:stepsize selection}, it is natural to derive a one-step descent inequality on the Lyapunov function defined in~\ref{eq:merit function:line-search-free}. Replacing \begin{equation}
    \label{eq:zeta_k}
\zeta^k:=\frac{\sqrt{(L^k)^2+2\sigma^k/c_1}-L^k}{2}>0,
\end{equation}
in~\eqref{eq:Lyapunov dynamic} in Lemma~\ref{lemma:Lyapunov dynamic} yields the following descent lemma
    \begin{lemma}[Descent Lemma]
\label{lemma:descent lemma:line-search-free}
In state of Theorem~\ref{thm:convergence of iterates}, then for any pair of saddle point $(\X^*,\Y^*)$ of Lagrangian function~\eqref{eq:lagarangian line-search-free}, it holds that for any $k\geq 0$ that,
\begin{equation}
    \label{eq:descent:line-search-free}
 V_*^{k+1}\leq V_*^k-\mathcal{D}_*^k,
\end{equation}
where 
\begin{equation}
    \label{eq:D^k}
\begin{aligned}
&\mathcal{D}_*^k
        =2(1-c_2)\gamma^k\alpha^k\Delta  \mathcal{L_*}^k
        \\&+\frac{1-c_1}{\sigma^k}\|\Y^{k}-\Y^{k+1}\|^2+\frac{1-2\alpha^kL^k}{2}\|\X^k-\X^{k-1}\|^2
        \\
        &+\frac{\alpha^k(\sqrt{(L^k)^2+2\sigma^k/c_1}-L^k)}{2}\|\X^k-\X^{k+1}\|_W^2\geq 0.
   \end{aligned}
\end{equation}
\end{lemma}

From Lemma~\ref{lemma:descent lemma:line-search-free}, we are able to show further that the sequence of iterates are within a compact set.
\begin{lemma}[Bounded Sequence]
    \label{lemma:bounded sequence}
In state of Theorem~\ref{thm:convergence of iterates}, then for any saddle points $(\X^*,\Y^*)$ of the Lagrangian function~\eqref{eq:lagarangian line-search-free}, it holds for any $k\geq 0$ that
    \begin{equation}
        \label{eq:bounded iterates}
        \|\X^k-\X^*\|^2+\frac{1}{\overline{\sigma}}\|\Y^k-\Y^*\|^2\leq V^0_*<\infty,
    \end{equation}
        and hence both $\{(\X^k,\D^k)\}_{k\geq 0}$ and $\{\Y^k\}_{k\geq 0}$ are bounded. Therefore, there exists $\widetilde{L}\in(0,\infty)$ such that $F$ is $\widetilde{L}-$smooth over the convex hull of $\cup_{t=0}^{\infty}[\X^{t-1},\X^t]$, meaning that $L^k\leq \widetilde{L}$ for any $k\geq 0$. Moreover, it can be guaranteed that
    \begin{equation}
        \label{eq:consecutive convergence}
     \begin{aligned}
\lim_{k\rightarrow\infty}\alpha^k(\X^k-\X^{k+1})=& \lim_{k\rightarrow\infty}\Y^k-\Y^{k+1}
=0.
     \end{aligned}
    \end{equation}
\end{lemma}

Observing~\eqref{eq:consecutive convergence}, it will be essential to guarantee that $\{\alpha^k\}$ is bounded from below to prove the convergence of iterates. The following Lemma~\ref{lemma:bounded stepsize related} verifies such guarantee.
\begin{lemma}[Bounded stepsize]
    \label{lemma:bounded stepsize related}
    Let $\{\alpha^k\}_{k\geq 0}$, $\{\gamma^k\}_{k\geq 0}$ be the sequences of stepsizes and ratio of two consecutive stepsizes respectively generated by Algorithm~\ref{alg:line-search-free} in state of Theorem~\ref{thm:convergence of iterates}. Then it holds that 
    \begin{equation}
        \label{eq: better stepsize bounded away from zero}
        \alpha^k\geq\underline{\alpha}:=\min\left\{\alpha^0,\left(\sqrt{\widetilde{L}^2+2\overline{\sigma}/c_1}+\widetilde{L}\right)^{-1}\right\};
    \end{equation}
    \begin{equation}
        \label{eq:bounded gamma}
       0<\gamma^k\leq \overline{\gamma}:=\frac{\sqrt{c_2^2+4}+c_2}{2}.
    \end{equation}
\end{lemma}

\begin{proof}
 ~\eqref{eq: better stepsize bounded away from zero} can be proved by induction. Note that by $\widetilde{L}-$smoothness of $F$ over the convex hull $\sup_{t=0}^{\infty}[\X^{t-1},\X^t]$, it holds that
 $L^k\leq \widetilde{L}$ for any $k\geq 0$ and thus by Asuumption~\ref{ass:convex para} \begin{equation}
     \label{eq:lower bound 1}
\frac{1}{\sqrt{(L^k)^2+2\sigma^k/c_1}+L^k}\geq \frac{1}{\sqrt{\widetilde{L}^2+2\overline{\sigma}/c_1}+\widetilde{L}}
 \end{equation}
 For $k=0$, $\alpha^0\geq\underline{\alpha}$ obviously. Assume~\eqref{eq: better stepsize bounded away from zero} holds for $k=K$, then at iteration $k=K+1$, we have 
 $$\alpha^{K+1}\overset{\eqref{eq:stepsize selection}\eqref{eq:lower bound 1}}{\geq}\min\left\{\frac{1}{\sqrt{\widetilde{L}^2+2\overline{\sigma}/c_1}+\widetilde{L}},\alpha^K\right\}\geq \underline{\alpha}.$$
 Therefore, by induction,~\eqref{eq: better stepsize bounded away from zero} holds for any $k\geq 0$. Now consider~\eqref{eq:bounded gamma}, note that for any $k\geq 0$, it holds that
 $$\begin{aligned}
     0<\gamma^k=&\frac{\alpha^k}{\alpha^{k-1}}\overset{\eqref{eq:stepsize selection}}{\leq}\frac{\sqrt{1+c_2\gamma^{k-1}}\alpha^{k-1}}{\alpha^{k-1}}=\sqrt{1+c_2\gamma^{k-1}}\\
     &\leq \sqrt{1+c_2\sqrt{1+c_2\gamma^{k-2}}}\leq\cdots\leq \frac{\sqrt{c_2^2+4}+c_2}{2}.
 \end{aligned}$$
\end{proof}

Finally, with Lemma~\ref{lemma:descent lemma:line-search-free},~\ref{lemma:bounded sequence} and~\ref{lemma:bounded stepsize related}, Theorem~\ref{thm:convergence of iterates} can be proved.
\begin{proof}[proof of Theorem 1]
    By Lemma~\ref{lemma:bounded sequence}, there exists a subsequence $\{N_k\}\subseteq\mathbb{N}_+$ such that $\{(\X^{N_k},\Y^{N_k},\D^{N_k})\}$ converges. We assume the cluster points are $\X^{\infty},\Y^{\infty}$ and $\PL\Y^{\infty}$ respectively. Combine Lemma~\ref{lemma:bounded sequence} and Lemma~\ref{lemma:bounded stepsize related} and by update~\eqref{eq:dual update}, we have
   \begin{equation}
       \label{eq:lemma 2.1 KKT}
\begin{aligned}
\nabla F(\X^{\infty})+\PL\Y^{\infty}=0.
\end{aligned}
   \end{equation}
   Note that $0<\sigma^0\leq\sigma^k\leq\overline{\sigma}$ and by Lemma~\ref{lemma:bounded stepsize related}, $\gamma^k\leq \overline{\gamma}$ for any $k\geq 0$, meaning $\sigma^k \gamma^k\alpha^k(\X^k-\X^{k-1})\rightarrow0$ and by~\eqref{eq:dual update},
   \begin{equation}
   \label{eq:consensus}
\lim_{k\rightarrow\infty}\alpha^{N_k}(I-W)\X^{N_k}=0\Rightarrow \PL\X^{\infty}=0.      
   \end{equation}
   Combining~\eqref{eq:lemma 2.1 KKT} and~\eqref{eq:consensus}, we conclude $(\X^{\infty},\Y^{\infty})$ is a pair of saddle points of~\eqref{eq:lagarangian line-search-free}. Therefore, by Lemma~\ref{lemma:descent lemma:line-search-free}, for any $k\geq 0$, we have
   $V_{\infty}^{k+1}\leq V_{\infty}^{k}.$
   Since $(\X^{\infty},\Y^{\infty})$ are cluster points, given $\varepsilon>0$, there exists $K>0$ such that $V_{\infty}^k<\varepsilon$ for any $k\geq K$. This completes the proof of~\eqref{eq:convergence of iterates}

 Apply Young's inequality on~\eqref{eq:auxiliary update} with $\varepsilon_k\in(0,1)$,
    \begin{equation}
        \label{eq:lemma2.4 inequality 1}
        \begin{aligned}
            &\|\Y^{k+1}-\Y^k\|^2\geq\left(1-\varepsilon_k\right)(1-c_1)\sigma^k(\alpha^k)^2\|\PL\X^k\|^2\\&-\left(\frac{1}{\varepsilon_k}-1\right)(1-c_1)\sigma^k(\gamma^k)^2(\alpha^k)^2\|\PL(\X^k-\X^{k-1})\|^2.
        \end{aligned}
    \end{equation}
Taking $\varepsilon_k:=\frac{(\alpha^k)^2}{M_1+(\alpha^k)^2}$ in~\eqref{eq:lemma2.4 inequality 1}, it holds
for any $k\geq0$,
    \begin{equation}
        \label{eq:PLX}
        \begin{aligned}
            &\frac{1-c_1}{\sigma^k}\|\Y^{k+1}-\Y^k\|^2
            \geq -\left(\frac{1}{2}-
        \alpha^kL^k\right)\|\X^k-\X^{k-1}\|^2\\
        &+(1-c_1)\sigma^0\frac{M_1}{M_1/(\underline{\alpha})^2+1}\|\PL\X^k\|^2.\\
        \end{aligned}
    \end{equation}
    Replace~\eqref{eq:PLX} in~\eqref{eq:descent:line-search-free} and by Assumption~\ref{ass:convex para}, for $k\geq 0$,
\begin{equation}
\label{eq:telescope}
   \begin{aligned}        V^k_\star-V^{k+1}_\star
        \geq M_2\left(\gamma^k\Delta  \mathcal{L}_\star^k+\gamma^k\|\PL\X^k\|^2\right),
   \end{aligned}
\end{equation}

   \noindent Note that $\mathcal{M}_\star(\X)$ is a convex function and $\mathcal{M}_{\star}(\X^k)=\Delta\mathcal{L}^k_\star+\|\PL X^k\|^2$. For any $k\geq 1$, the convexity of $\mathcal{M}_\star(\X)$ yields 
   \begin{equation}
       \begin{aligned}
           \mathcal{M}_\star(\overline{\X}^k)
           \overset{\eqref{eq:telescope}}{\leq} \frac{1}{\theta^k}\sum_{t=0}^{k-1}\frac{V^t_\star-V^{t+1}_\star}{M_2}
           \leq \frac{V^0_\star/M_2}{\theta^k}.
       \end{aligned}
   \end{equation}
  By the definition of $\gamma^k$ in~\eqref{eq:stepsize selection}, for any $k\geq 2$,
  $$\frac{\theta^k}{k}\geq\left(\prod_{t=0}^{k-1}\gamma^t\right)^{1/k}=\left(\gamma^0\prod_{t=1}^{k-1}\frac{\alpha^t}{\alpha^{t-1}}\right)^{1/k}=\left(\frac{\alpha^{k-1}}{\alpha^0}\right)^{1/k},$$
   By~\eqref{eq: better stepsize bounded away from zero}, we have
  $\theta^k\geq k\cdot\left(\frac{\underline{\alpha}}{\alpha^0}\right)^{1/k}$. 
Since $\underline{\alpha}/\alpha^0\leq 1$, there exists $K\geq 2$ such that $(\underline{\alpha}/\alpha^0)^{1/k}\geq 1/2$, $\forall k\geq K$ and hence $\frac{1}{\theta^k}\leq \mathcal{O}(1/k)$, which completes the proof. 
\end{proof}
\section{Proof of Theorem~\ref{thm:contraction}}
\label{appx:proof of theorem 2}
\begin{proof}
The convergence of iterates can be proved similar to that in Theorem~\ref{thm:contraction} and we focus on the linear convergence here. Define a new merit function: for $k\geq 0$,
\begin{equation}
    \label{eq:merit linear}\begin{aligned}
&\widetilde{V}^{k+1}=\|\X^{k+1}-\X^*\|^2+\left(\frac{1}{\sigma^k}+\vartheta^{k}\right)\|\Y^{k+1}-\Y^*\|^2\\
&+(2\gamma^{k+1}\alpha^{k+1}+2(1-c_2)\gamma^k\alpha^k)\Delta  \mathcal{L}_*^k+\frac{1}{2}\|\X^{k+1}-\X^k\|^2
,
\end{aligned}
\end{equation}
where\begin{equation}
    \label{eq:vartheta}
\vartheta^k:=(1-\lambda_2(W))\cdot\min\left\{\frac{\sigma^k(\alpha^k)^4\lambda_m(W)}{2c_1},\frac{\alpha^k}{4\widetilde{L}}\right\}.
\end{equation}
    Redo the analysis with $\mu>0$ for Lemma~\ref{lemma:descent lemma:line-search-free} with $\mu>0$. Then by Lemma~\ref{lemma:bounded stepsize related}, it holds that for any $k\geq 0$
\begin{equation}
    \label{eq:enhanced descent}
    \begin{aligned}
        V_*^{k+1}\leq &V_*^k-\frac{\alpha^k}{2\widetilde{L}}\|\nabla F(\X^k)-\nabla F(\X^*)\|^2\\
        &-\frac{\alpha^k\widetilde{\mu}}{2}\|\X^k-\X^*\|^2-\mathcal{D}_*^k.
    \end{aligned}
\end{equation}
Note that by stepsize selection rule~\eqref{eq:stepsize selection},
\begin{equation}
    \label{eq:thm2.3 inequality 1}
    \frac{\alpha^k}{2}\left(\sqrt{(L^k)^2+2\sigma^k/c_1}-L^k\right)
    \geq \frac{(\alpha^k)^2\sigma^k}{c_1}.
\end{equation}
By update of primal variable~\eqref{eq:primal update}, the fact that $\PL\Y^*+\nabla F(\X^*)=0$ and $\Y^k\in\texttt{span}(\PL)$, we have for any $k\geq 0$,
\begin{equation}
         \label{eq:dual gap bound}
       \begin{aligned}
             &\|\Y^{k+1}-\Y^*\|^2\leq \frac{2}{1-\lambda_2(W)}\|\nabla F(\X^k)-\nabla F(\X^*)\|^2\\
             &+\frac{2}{(\alpha^k)^2(1-\lambda_2(W))\lambda_m(W)}\|\X^k-\X^{k+1}\|_W^2.
       \end{aligned}
     \end{equation}
Then by combining~\eqref{eq:enhanced descent},~\eqref{eq:thm2.3 inequality 1} and~\eqref{eq:dual gap bound}, for any $k\geq 1$, it holds that
  \begin{equation}
      \label{eq:contraction}
      \widetilde{V}^{k+1}\leq (1-\rho_k)\widetilde{V}^k,
  \end{equation}
  where
  \begin{equation}
      \label{eq:rho_k}
      \begin{aligned}
          \rho_k:=&\min\left\{\frac{\alpha^k\widetilde{\mu}}{2},1-\frac{1}{\sigma^k/\sigma^{k-1}+\sigma^k\vartheta^{k-1}},1-2\alpha^kL^k,\right.\\
          &\left.1-\frac{2\gamma^k\alpha^k}{2\gamma^{k}\alpha^{k}+2(1-c_2)\gamma^{k-1}\alpha^{k-1}}\right\}<1,
      \end{aligned}
  \end{equation}
  and $\vartheta^k$ is defined in~\eqref{eq:vartheta}.   
   Note that under Assumption~\ref{ass:strong convx para} we have $\sigma^k=\sigma/(\alpha^k)^2$ and $\alpha^k\leq \widetilde{\gamma}_k\alpha^{k-1}$ with
   $$\widetilde{\gamma}^k:=((k+\beta_1)/(k+1))^{\beta_2}.$$
   By~\eqref{eq:stepsize selection linear implementation}, for any $k\geq 1$,
    \begin{subequations}
    \label{eq:cor3.1 inequalities}
        \begin{align}
        \label{eq:cor3.1 inequality 1}
         \frac{\alpha^k\widetilde{\mu}}{2} \geq\widetilde{\rho}_1; \\
         \label{eq:cor3.1 inequality 2}
         1-\frac{1}{\sigma^k/\sigma^{k-1}+\sigma^k\vartheta^{k-1}}&\geq1-(\gamma^k)^2\widetilde{\rho_2};\\
         \label{eq:cor2.3 inequality 3}
         1-2\alpha^kL^k\geq\widetilde{\rho}_3:=\frac{2\sigma}{c_1}<1; \\
         \label{eq:cor2.3 inequality 4}
         1-\frac{2\gamma^k\alpha^k}{2\gamma^k\alpha^k+2(1-c_2)\gamma^{k-1}\alpha^{k-1}}&\geq\widetilde{\rho}_4,
        \end{align}
    \end{subequations}
    where $$\begin{aligned} 
    0<\widetilde{\rho}_1&:= \max\left\{1-\frac{\alpha_0\widetilde{\mu}}{2},1-\frac{1/2-\sigma/c_1}{2\widetilde{\kappa}}\right\}\\
    0<\widetilde{\rho}_2:=&\min\left\{\frac{1}{1+\sigma^2(1-\lambda_2(W))\lambda_m(W)/(2c_1)},\right.\\
    &\left.\frac{1}{1+\sigma(1-\lambda_2(W))/((2-2\sigma/c_1)\widetilde{\kappa})}\right\},\\
   0< \widetilde{\rho}_4:=&1-\frac{1}{1+(1-c_2)/\left(\overline{\gamma}^2\max\{\widetilde{\kappa},(1/2-\sigma/c_1)/(\alpha^0\mu)\}\right)}.
    \end{aligned}$$ 
    
    Then by the fact that $\widetilde{\gamma}^k\geq\max\{1,\gamma^k\}$ for any $k\geq 1$, taking $\rho:=\min\{\widetilde{\rho}_1,1-\widetilde{\rho}_2,\widetilde{\rho}_3,\widetilde{\rho}_4\}$ yields: any $k\geq 1$,
    \begin{equation}
    \label{eq:rho bound}
        \begin{aligned}
        1-\rho_k\leq& 1-\min\left\{\widetilde{\rho}_1,1-(\gamma^k)^2\widetilde{\rho}_2,\widetilde{\rho}_3,\widetilde{\rho_4}\right\}\\
        =&\max\{1-\widetilde{\rho}_1,(\gamma^k)^2\widetilde{\rho}_2,1-\widetilde{\rho}_3,1-\widetilde{\rho}_4\}\\
        \leq& (\widetilde{\gamma}^k)^2\max\{1-\widetilde{\rho}_1,\widetilde{\rho}_2,1-\widetilde{\rho_3},1-\widetilde{\rho}_4\}\\
        =&(\widetilde{\gamma}^k)^2(1-\rho)
    \end{aligned}.
    \end{equation} 
    
    By Theorem~\ref{thm:contraction}, we have that for $k\geq 1$
    \begin{equation}
        \label{eq:new contraction}
        \widetilde{V}^{k+1}\leq (1-\rho_k)\widetilde{V}^k\overset{\eqref{eq:rho bound}}{\leq} (\widetilde{\gamma}^k)^2(1-\rho) \widetilde{V}^k.
    \end{equation}
Therefore by Assumption~\ref{ass:strong convx para}, there exists $K'$ such that for $N_{\varepsilon}>K'$, it holds that
$$\|\X^{N_{\varepsilon}+1}-\X^*\|^2\leq \left(\prod_{k=1}^{N_{\varepsilon}+1}\widetilde{\gamma}^k\right)^2(1-\rho)^{N_{\varepsilon}}\widetilde{V}^1\leq (1-\rho/2)^{N_{\varepsilon}}\widetilde{V}^1,$$
where
\begin{equation}
    \label{eq:rho}
    \rho=\Omega\left(\frac{1}{1+\widetilde{\kappa}/(1-\lambda_2(W))}\right).
\end{equation}
\end{proof}

\section{Proofs for Section~\ref{sec:local-min}}
\label{appx:section 3}
\begin{proof}[proof of Lemma~\ref{lemma:conditioned stepsize consistency}]
Since $\mathcal{K}$ is finite, there exists $K_1>0$ such that $\hat{\alpha}_i^k>\pi^k(\alpha_i^k)$ for any $i\in[m]$. Then by~\eqref{eq:tilde alpha 1}, for $k\geq K_1$, it holds for any $i\in[m]$ that $$\widetilde{\alpha}_i^k=\min\{\sqrt{1+c_2\gamma_i^{k-1}}\alpha_i^{k-1},\pi^k(\alpha_i^{k-1})\}\geq \alpha_i^{k-1}.$$
With $\gamma_i^k$ set in~\eqref{eq:stepsize selection local 1}, it holds that $\gamma_i^k\geq 1$ for any $k\geq K_1+1$. Note that the sequence $\{\X^k\}$ is bounded, meaning that $\hat{\alpha}_i^k$ is bounded from below, it can be easily shown by induction that $\{\alpha_i^k\}$ selected by~\eqref{eq:stepsize selection local 1} is bounded from below. Then by Assumption~\ref{ass:Pi local},there exists $K_2\geq K_1$ such that for any $k\geq K_2$, it holds for any $i\in[m]$ that
$$\widetilde{\alpha}_i^k=\min\{\sqrt{1+c_2\gamma_i^{k-1}}\alpha_i^{k-1},\pi^k(\alpha_i^{k-1})\}=\pi^k(\alpha_i^{k-1}).$$
Therefore, the stepsizes across the network will converge in $d_{\mathcal{G}}$ number of iterations where $d_{\mathcal{G}}$ denotes the diameter of the network. After that, $\widetilde{\alpha}_i^k=\alpha_i^k$ and hence $\gamma_i^k=\alpha_i^k/\alpha_i^{k-1}$. Take $K=K_2+d_{\mathcal{G}}$ and we complete the proof.
\end{proof}

\begin{proof}[proof of Lemma~\ref{lemma:finite K}]
Since $\{\X^k\}$ is bounded, $\hat{\alpha}_i^k\geq \underline{\hat{\alpha}}$ for some $\underline{\hat{\alpha}}$ for any $i\in[m]$ and $k\geq 0$. Since $\{\delta^k\}$ is summable, there exists $K_3>0$ such that
$\sum_{t=k}^\infty\delta^t<\frac{\underline{\hat{\alpha}}}{2}$ for any $k\geq K_3$. Suppose on the contrary that there exists $i\in[m]$ that $\{\alpha_i^k\}$ has a subsequence that converge to zero, then there exists $k_1\geq K_3$ such that $\alpha_i^{k_1}\leq \frac{\hat{\alpha}}{2}$.  By~\eqref{eq:stepsize selection local 1}, it holds that for any $k>k_1+d_{\mathcal{G}}$
$$\pi^{k}(\alpha_j^{k-1})\leq \alpha_i^{k_1}+\sum_{t=k_1}^{k}\delta^t< \underline{\hat{\alpha}}\leq \hat{\alpha}_j^{k},\quad \forall j\in[m].$$
Therefore, it holds for any $k\geq k_1+d_{\mathcal{G}}$ and $\forall j\in[m]$ that
$$\widetilde{\alpha}_j^k=\min\left\{\pi^{k}(\alpha_j^{k-1}),\sqrt{1+c_2\gamma^{k-1}_i}\alpha_j^{k-1}\right\}\geq \alpha_j^{k-1}.$$
Meaning that for any $k\geq k_1+d_{\mathcal{G}}$, we have $$\min_{j\in[m]}\alpha_j^k=\min_{j\in[m]}\widetilde{\alpha}_j^k\geq \min_{j\in[m]}\alpha_{j}^{k-1}$$ and thus the sequence $\{\min_{j\in[m]}\alpha_j^k\}$ is bounded away from zero, contradict with our assumption that there exists $i\in[m]$ where ${\alpha_i^k}$ has a subsequence converging to zero. Therefore, the sequence $\{\alpha_i^k\}$ is bounded away from zero for any $i\in[m]$. We write $\alpha_i^k\geq\underline{\alpha}$ for any $i\in[m]$ and $k\geq 0$.

It remains to prove that $|\mathcal{K}|<\infty$. Still, assume on the contrary that $|\mathcal{K}|=\infty$. There there exists at least an $i\in[m]$ such that the subset $\mathcal{K}_i:=\{k:\hat{\alpha}_i^k\leq \pi^k(\alpha_i^{k-1})\}$ is infinite.
By the summability of $\{\delta^k\}$, there exists $K_4\geq 0$ such that $\sum_{t=k}^{\infty}\delta^t\leq \frac{1-\eta}{2}\underline{\alpha}$. Then for any $k_4>k_3>K_4$ and $k_3,k_4\in \mathcal{K}$, it holds by~\eqref{eq:stepsize selection local 1} that
$$\alpha_i^{k_4}\leq \eta\alpha_i^{k_3}+\sum_{t=k_3+1}^{k_4}\delta^t\leq \frac{1+\eta}{2}\alpha_i^{k_3}<\alpha_i^{k_3},$$
meaning that there exists subsequence of $\{\alpha_i^k\}$ converging to zero, contradiction and hence we have $|\mathcal{K}|<\infty$.
\end{proof}

Combining Lemma~\ref{lemma:conditioned stepsize consistency} and Lemma~\ref{lemma:finite K}, we can conclude that all the stepsizes over the network will become homogeneous after a finite number of iterations, as we stated in Lemma~\ref{lemma:stepsize stablization}. 

\begin{lemma}
\label{lemma:stepsize stablization}
    Let $\{\Lambda^k\}_{k\geq 0}$, $\{\Gamma^k\}_{k\geq 0}$ be the sequences of stepsizes and ratio of two consecutive stepsizes generated by Algorithm~\ref{alg:line-search-free_local}. If $\{\pi_k\}_{k\geq 0}$ and $\eta\in(0,1)$ are chosen such that under Assumption~\ref{ass:Pi local}. Then if the sequence $\{\alpha^k_i\}$ is lower bounded, there exists $K>0$ and for any $k\geq K$, it holds for some $\alpha^k,\gamma^k>0$ that
    \begin{equation}
        \label{eq:the same stepsize}
        \Lambda^k=\alpha^k I,\quad \Gamma^k=\gamma^k I,\quad \gamma^k=\frac{\alpha^k}{\alpha^{k-1}},
    \end{equation}
    where $\alpha^k$ satifies~\eqref{eq:stepsize selection} with $\sigma^k:=\min_{i\in[m]}\sigma_i^k$.
\end{lemma}
Then after such finite number of iterations, Theorem~\ref{thm:iterates convergence_local} and~\ref{thm:contraction local} can be proved naturally.

\end{document}